\documentclass{article}

\usepackage{arxiv}
\usepackage{fancyhdr}
\usepackage[utf8]{inputenc} 
\usepackage[T1]{fontenc}    
\usepackage{hyperref}       
\usepackage{url}            
\usepackage{booktabs}       
\usepackage{amsfonts}       
\usepackage{nicefrac}       
\usepackage{microtype}      
\usepackage{lipsum}
\usepackage{amsthm}
\usepackage{graphicx}
\usepackage{natbib}
\usepackage{color}
\usepackage{doi}
\usepackage{mathtools}
\usepackage{amssymb}
\usepackage{bm}
\usepackage{float}
\usepackage{amsmath}
\usepackage{empheq}
\usepackage{subcaption}
\usepackage{tabularx,booktabs,array,adjustbox}
\usepackage{array}
\usepackage{cancel}
\usepackage{placeins}

\renewcommand{\headeright}{\textit{arXiv Preprint}}

\renewcommand{\shorttitle}{}

\renewcommand{\undertitle}{}

\title{Space-time discretization for barotropic flow stemming from a multisymplectic  variational formulation}

\author{
    {Mukthesh Mahadev}\\
	Department of Flow Physics \& Technology \\
	Delft University of Technology\\
	Delft, Netherlands \\
	\texttt{mukthesh.mahadev@outlook.com} \\
	\And
	\href{https://orcid.org/0000-0003-2539-642X}{\includegraphics[scale=0.06]{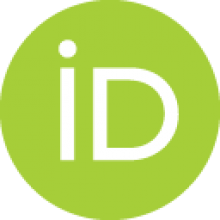}\hspace{1mm}Marc Gerritsma} \\
	Department of Flow Physics \& Technology \\
	Delft University of Technology\\
	Delft, Netherlands \\
	\texttt{M.I.Gerritsma@tudelft.nl} \\
}

\begin{document}

\maketitle

\begin{abstract}
This study proposes and analyses a novel higher-order, structure-preserving discretization method for inviscid barotropic flows from a Lagrangian perspective. The method is built on a multisymplectic variational principle discretized over a full space-time domain. Flow variables are encoded on a staggered space-time mesh, leveraging the principles of mimetic spectral element discretization. Unlike standard Lagrangian methods, which are prone to mesh distortion, this framework computes fluid deformations in a fixed reference configuration and systematically maps them to the physical domain via the Piola-Kirchhoff stress. Further, the structure-preserving design ensures that the discrete analogues of the fundamental conservation laws for mass, momentum, and energy are satisfied up to machine precision. The formulation also inherently handles low-Mach number flows without specialized preconditioning. Numerical experiments on expansion and compression flows confirm the accuracy, stability, and exact conservation properties of the discretization.
\end{abstract}


\section{Introduction}

In computational fluid dynamics, developing numerical methods that preserve crucial conservation properties over long time integration remains a significant challenge. Many approaches may match experimental data for specific cases but fail to conserve invariants such as energy, momentum, or mass, leading to numerical drift and/or instability. This has driven interest in structure-preserving frameworks that capture the geometric and variational foundations of fluid behaviour in space and time. Inviscid barotropic flow, where the pressure is dependent solely on the mass density, is one such class of flows governed by variational principles. A Lagrangian viewpoint, which tracks individual fluid parcels, is natural for these flows. However, traditional Lagrangian methods that deform a computational mesh in physical space are notoriously susceptible to severe mesh distortion and numerical instability. To avoid this, we employ a formulation defined in a fixed reference configuration, where fluid motion is captured through a flow map and dynamics is governed by the First Piola-Kirchhoff stress tensor ($\widetilde{P}$) rather than the Cauchy stress in the deformed domain \cite{Reddy}.

This work adopts the framework of mimetic discretization and Finite Element Exterior Calculus (FEEC) \cite{feec, bochev2006principles, hiptmair2002finite}, which is built on direct discretization of the underlying geometric structures of physical laws \cite{tonti2013mathematical, TONTI20141260}. The degrees of freedom for the kinematic variables are represented as integral values in the nodes, edges, faces and volumes of the mesh, depending on the physical character of the variable. We will refer to these variables as the primal variables. The dynamic variables are represented by dual fields that act on the primary variables. Continuous differential operators (grad, curl, div) will be represented by discrete topological incidence matrices, $\mathbb{E}$, that satisfy the de Rham complex at the discrete level. The constitutive laws relate primal to dual variables, i.e kinematic variables to dynamic variables. These relations will contain weighted mass matrices, $\mathbb{M}_k$, where $k$ denotes some material parameter, \cite{bochev2006principles,MimeticFramework2011}. This structure-preserving paradigm has been successfully applied to various problems as seen in \cite{palha_thesis,palha2017mass,zhang2022mimetic}. In this work, we extend this framework to a full space-time discretization of Lagrangian barotropic flow.

The primary contribution of this work is to develop a  multisymplectic space-time mimetic method with,
\begin{itemize}
    \item A full space-time (2+1D) discretization of the Lagrangian barotropic flow equations.
    \item A formulation based on the first Piola-Kirchhoff stress tensor that evaluates physical stress components by evaluating them in the reference domain.
    \item A discrete weak formulation that is shown to \textit{exactly} conserve mass, linear momentum, angular momentum, and total energy.
    \item Numerical verification of these conservation properties and a demonstration of the method's stability for low-Mach number expansion and compression flows.
\end{itemize}

The outline of this paper is as follows: In Section~\ref{sec:Largrangina_formulation} the weak formulation is derived from a Lagrangian density. In Section~\ref{sec:MimeticDiscretization} basis functions and incidence matrices are introduced. In Section~\ref{sec:SpaceTimeFormulation} the representation of all physical variables in the space-time grid is explained.  Discretization of the governing equations is discussed in Section~\ref{sec:algebraic_reduction}. Derivation of the conservation properties for the discrete system is addressed in Section~\ref{sec:conservation_properties}. The corroboration of these conservation properties is given in Section~\ref{sec:NumericalResults}. Concluding remarks can be found in Section~\ref{sec:Conclusions}.

\section{Lagrangian Formulation of Barotropic Flow}
\label{sec:Largrangina_formulation}

\begin{figure}[H]
    \centering
    \includegraphics[width=0.4\linewidth]{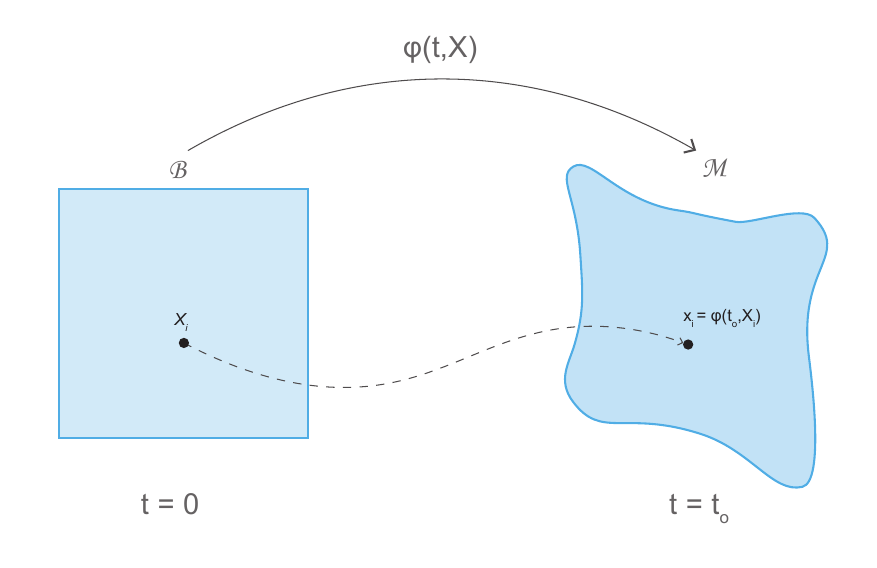}
    \caption{Reference and Deformed Configurations with Flow Map \(\varphi(t, X)\)}
    \label{fig:flowmap_defnition}
\end{figure}

To derive the barotropic flow equations, we first establish the kinematic and dynamic principles of the Lagrangian framework \cite{DemouresGayBalmaz2021,marsden1994mathematical,Reddy}. The reference configuration as shown in Figure \ref{fig:flowmap_defnition}, is denoted as the undeformed state of the fluid, is a compact domain \( \mathcal{B} \subset \mathbb{R}^n \), where each point \( \bm{X} \in \mathcal{B} \) is a unique particle label. The motion of the fluid is then described by the flow map, \( \bm{\varphi} : \mathcal{B} \times \mathbb{R}^+ \rightarrow \mathcal{M} \), which provides the physical position \( \bm{x} = \bm{\varphi}(t, \bm{X}) \) of particle \( \bm{X} \) at time \( t \). The set of all particle positions at time $t$ forms the deformed configuration $\mathcal{M}(t)$.

Local deformation is quantified by the deformation gradient tensor, $\bm{F} = \bm{\nabla}_X \bm{\varphi}$, with components $F_{ij} = \partial \varphi^i / \partial X_j$. The determinant of this tensor, $J = \det(\bm{F})$, relates the reference density $\rho_o$ to the deformed density $\rho$ through the conservation of mass,
\begin{equation}
    J = \frac{\rho_o}{\rho(t, \bm{X})}.
    \label{eq:Jacobian_determinant}
\end{equation}

A fluid is termed \textit{barotropic} when its pressure is a sole function of its mass density, $p = p(\rho)$. Consequently, the specific internal energy $W$ is also only a function of $\rho$, and thus of $J$. From the first law of thermodynamics, $dW = -p \, d(1/\rho)$, we can define a material pressure $P_W$ (work-conjugate to $J$) that represents the change in internal energy with respect to the change in volume ratio. This relationship is given by,
\begin{equation}
    P_W(\rho_o, J) = - \rho_o \frac{\partial W(\rho_o, J)}{\partial J}.
    \label{eq:material_pressure}
\end{equation}

\subsection{Hamilton's Principle and the Variational Formulation}
The equations of motion are derived from the well known Hamilton's principle of stationary action, $\delta \mathcal{A} = 0$. The action $\mathcal{A}$ is the integral of the Lagrangian density $\mathcal{L} = T - \mathcal{V}$, where \ $T = \frac{1}{2}\rho_o\left|\dot{\bm{\varphi}} \right|^2$ is the kinetic energy density and $\mathcal{V} = \rho_o W(\rho_o, J) + \rho_o \Pi(\bm{\varphi})$ \ is the potential energy density. The full Lagrangian density is therefore \cite{DemouresGayBalmaz2021}:
\begin{equation}
\mathcal{L}(\bm{\varphi},\dot{\bm{\varphi}},\nabla \bm{\varphi}) = \frac{1}{2}\rho_0 |\dot{\bm{\varphi}}|^2 - \rho_0 W(\nabla \bm{\varphi}) - \rho_0 \Pi(\bm{\varphi})
\label{eq:lagrangian}
\end{equation}
The stationary point of the action is found by taking the first variation of the action integral and then setting it to zero as stated by Hamilton's principle, this yields the well known Euler-Lagrange equations, which state the strong form of the momentum balance,
\begin{equation}
\frac{\partial}{\partial t}\left( \frac{\partial\mathcal{L}}{\partial\dot{\bm{\varphi}}} \right) + \frac{\partial}{\partial X_j}\left( \frac{\partial\mathcal{L}}{\partial\nabla_j\bm{\varphi}} \right) - \frac{\partial\mathcal{L}}{\partial{\bm{\varphi}}} = 0.
\label{eq:euler-lagrange}
\end{equation}
Substituting \eqref{eq:lagrangian} into \eqref{eq:euler-lagrange} yields the material acceleration $\rho_o \ddot{\bm{\varphi}}$ and the external force $\rho_o \frac{\partial \Pi}{\partial \bm{\varphi}}$ for the first and the last term respectively. The stress term is derived by applying the chain rule to the internal energy $W$,
\begin{gather*}
\frac{\partial\mathcal{L}}{\partial\nabla_j\varphi^i} = \frac{\partial}{\partial\nabla_j\varphi^i} \left(-\rho_o W(\rho_o,J)\right) = \left(-\rho_o\frac{\partial W}{\partial J}\right) \left(\frac{\partial J}{\partial\nabla_j\varphi^i}\right).
\end{gather*}
From \eqref{eq:material_pressure}, we identify the material pressure as $P_W = - \rho_o \frac{\partial W}{\partial J}$. The second term, the derivative of the determinant $J = \det(\bm{F})$, is given by Jacobi's formula,
\begin{equation*}
\frac{\partial J}{\partial\nabla_j\varphi^i} = \frac{\partial(\det (\bm{F}))}{\partial\bm{F}_{ij}} = \det(\bm{F}) (\bm{F}^{-1})_{ji} = J\bm{F}^{-T}.
\end{equation*}
Combining these, defines the first Piola-Kirchhoff stress tensor $\bm{P} = \frac{\partial\mathcal{L}}{\partial(\nabla \bm{\varphi})} =  P_W J \bm{F}^{-T} $. Substituting this back into \eqref{eq:euler-lagrange} gives the strong form of the Lagrangian barotropic flow equations,

\begin{equation}
\begin{aligned}
    \rho_o\ddot{\varphi}_i +\frac{\partial}{\partial X_j}\left( P_W J\bm{F}^{-T} \right)_{ij}+ \rho_o\frac{\partial\Pi}{\partial\varphi_i} = 0\;.
\end{aligned}
\label{eq:Lagrangian_Barotropic_Equations}
\end{equation}

Although \eqref{eq:Lagrangian_Barotropic_Equations} is the strong form, our numerical method is directly based on the first variation of the action $\delta\mathcal{A}$ to remain in an integral setting,

\begin{equation}
\delta \mathcal{A} = \int\limits_{0}^{T} \int\limits_{\mathcal{B}} \left[ \rho_0 \big ( \dot{\bm{\varphi}}, \delta\dot{\bm{\varphi}} \big ) - \rho_0 \frac{\partial W}{\partial (\nabla \bm{\varphi})} : \nabla(\delta\bm{\varphi}) - \rho_0 \left \langle \frac{\partial \Pi}{\partial \bm{\varphi}}, \delta\bm{\varphi} \right \rangle \right] \, \mathrm{d}\mathcal{B} \ \mathrm{d}t = 0\;.
\label{eq:weak_form}
\end{equation}

To create a first-order system suitable for a space-time discretization, we introduce the canonical momentum $\bm{\pi}$ as an independent field and enforce its definition, $\widetilde{\bm{\pi}} = \rho_0 \dot{\bm{\varphi}}$, with a second weak equation. Further, using the previously derived definition of the Piola-Kirchhoff stress tensor being equal to the derivative of internal energy with respect to the gradient of the flow map, the final two-field variational formulation can be written as,

\begin{empheq}[]{equation}
\begin{gathered}
\int\limits_{0}^{T}\!\int\limits_{\mathcal{B}}
\Big[
  \langle \widetilde{\bm{\pi}}, \delta\dot{\bm{\varphi}} \rangle
  + P_{W}J\bm{F}^{-T} : \nabla(\delta\bm{\varphi})
  - \rho_0 \langle \tfrac{\partial \Pi}{\partial \bm{\varphi}}, \delta\bm{\varphi} \rangle
\Big]\,
\mathrm{d}\mathcal{B} \ \mathrm{d}t
= 0 \quad \forall\, \delta\bm{\varphi} \\
\int\limits_{0}^{T}\!\int\limits_{\mathcal{B}}
\Big[
  \big ( \widetilde{\bm{\pi}},\delta\widetilde{\bm{\pi}} \big ) 
  -\rho_0 \left \langle \dot{\bm{\varphi}}, \delta\widetilde{\bm{\pi}} \right \rangle
\Big]\,
\mathrm{d}\mathcal{B} \ \mathrm{d}t
= 0 \quad \forall\, \delta\bm{\pi}
\end{gathered}
\label{eq:pre_final_weak_formulation}
\end{empheq}

The above equations, after being supplemented with the temporal (initial) and spatial boundary conditions, are the equations that are discretised in this work, that is, to find $\bm{\varphi} \ \text{and} \  \bm{\pi}$ such that they hold for all variations $(\delta\bm{\varphi}, \delta\bm{\pi})$. Further, $\langle \cdot , \cdot \rangle$ denotes duality pairing between the momentum and velocity variables, and $(\cdot,\cdot)$ and $\cdot : \cdot$ denote the vector and tensor inner products, respectively. It will turn out later in this paper, that duality pairing only involves the degrees of freedom and there is no dependence on the basis functions, while the inner product depends explicitly on the basis functions. For the remainder of this paper, we set the potential $\Pi(\bm{\varphi})$ to zero. So the final weak formulation becomes \\

\begingroup
\small 

\setlength{\fboxsep}{3pt}    
\setlength{\fboxrule}{0.5pt} 

\begin{empheq}[box=\fbox]{equation}
\begin{gathered}
\int\limits_{0}^{T}\!\int\limits_{\mathcal{B}}
\Big[
  \langle \widetilde{\bm{\pi}}, \delta\dot{\bm{\varphi}} \rangle
  + P_{W}J\bm{F}^{-T} : \nabla(\delta\bm{\varphi})
\Big]\,
\mathrm{d}\mathcal{B} \ \mathrm{d}t - \int_\mathcal{B} \left \langle \widetilde{\bm{\pi}}, \delta \bm{\varphi} \right \rangle \big |_{t=0}^{t=T} \mathrm{d} \mathcal{B}
= \int_0^T \int_{\partial \mathcal{B}} \left ( \widehat{P}_{ext} J \bm{F}^{-T} \bm{n}, \delta \bm{\varphi} \right ) \mathrm{d}(\partial \mathcal{B}) \mathrm{d}t \quad \forall\, \delta\bm{\varphi} \\
\int\limits_{0}^{T}\!\int\limits_{\mathcal{B}}
\Big[
  \big ( \widetilde{\bm{\pi}},\delta\widetilde{\bm{\pi}} \big )
  -\rho_0 \left \langle \dot{\bm{\varphi}}, \delta\widetilde{\bm{\pi}} \right \rangle
\Big]\,
\mathrm{d}\mathcal{B} \ \mathrm{d}t
= 0 \quad \forall\, \delta\bm{\pi}
\end{gathered}
\label{eq:final_weak_formulation}
\end{empheq}

\endgroup

\section{Mimetic Discretization}\label{sec:MimeticDiscretization}

To numerically solve the variational formulation derived in \eqref{eq:final_weak_formulation}, we must construct a discrete finite dimensional subspace that preserves the geometric structure of the continuous equations. To achieve this, we employ the Mimetic Spectral Element Method \cite{MimeticFramework2011,introCompatiblemethods}. In this framework, variables are not simply defined at grid points, rather they are associated with geometric entities of the mesh nodes, edges, faces, or volumes based on their physical interpretation. Consequently, we cannot rely solely on standard nodal basis functions to represent these fields. Instead, we must construct a hierarchy of compatible basis functions that correspond to these specific geometric entities. This section introduces the construction of these nodal and edge-based polynomial spaces.

\paragraph{Nodal (Lagrange) Basis.}
Scalar fields are approximated using the nodal Lagrange basis $\{h_i(\xi)\}$ associated with the GLL nodes $\{\xi_i\}$ on the reference interval $[-1,1]$. These basis functions satisfy the Kronecker delta property $h_i(\xi_j) = \delta_{ij}$. The finite dimensional representation $\phi_h(\xi)$ from the degrees of freedom $\overline{\phi}$ (which stores nodal values $\phi_i=\phi(\xi_i)$) is:
\begin{equation}
  \label{eq:0form_reconstruction}
  \phi_h(\xi)
  \;=\;
  \sum_{i=0}^N
  \phi_i \, h_i(\xi) \; .
\end{equation}

\paragraph{Edge Basis.}
Alternatively, functions can be represented by edge basis functions $\{e_i(\xi)\}$ \cite{gerritsma2011edge}. These are constructed to satisfy $\int_{\xi_{j-1}}^{\xi_j} e_i(\xi) \, \mathrm{d}\xi = \delta_{ij}$, which means that the degrees of freedom $\overline{\psi}$ represent the integral values, $\psi_i = \int \limits_{\xi_{i-1}}^{\xi_i} \psi(\xi) \ \mathrm{d}\xi $ . The finite dimensional representation for $\psi_h(\xi)$ is,
\begin{equation}
  \label{eq:1form_reconstruction}
  \psi_h(\xi)
  \;=\;
  \sum_{i=1}^N
  \psi_i \, e_i(\xi).
\end{equation}

\begin{figure}[h!]
    \centering
    \begin{subfigure}[b]{0.35\textwidth}
        \centering
        \includegraphics[width=\linewidth]{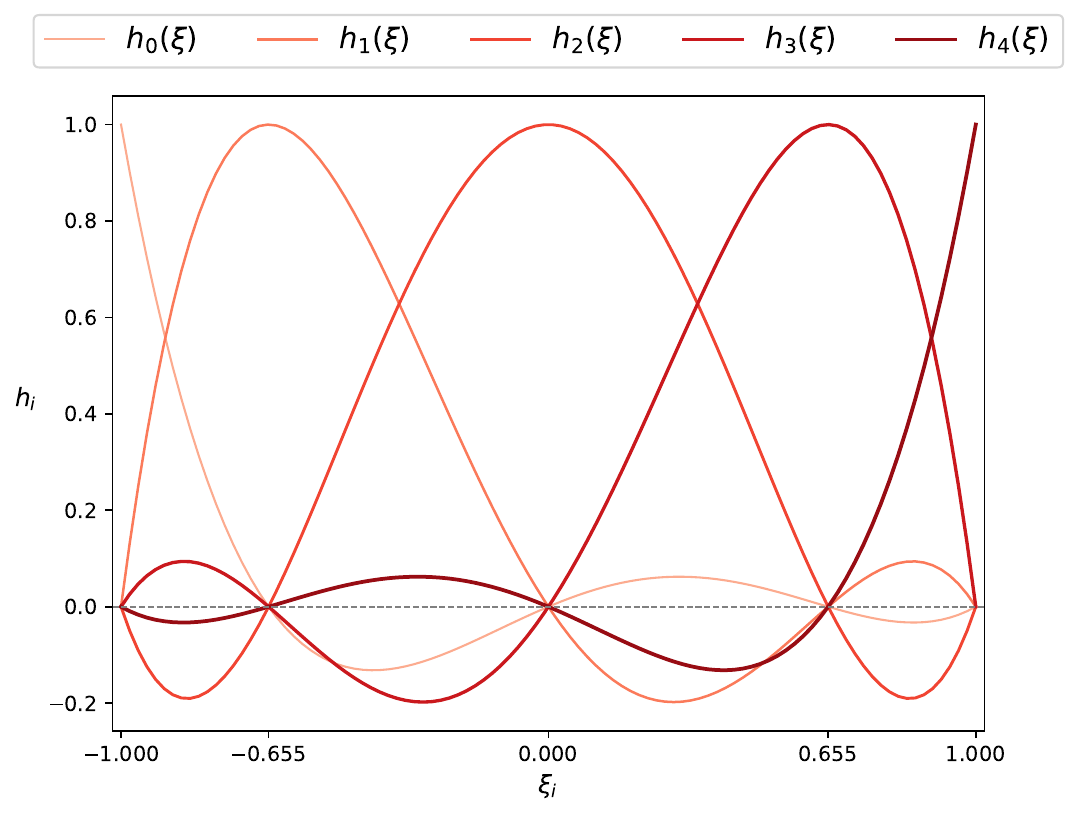}
        \caption{Nodal Lagrange polynomial basis}
    \end{subfigure}
    \hspace{2cm}
    \begin{subfigure}[b]{0.32\textwidth}
        \centering
        \includegraphics[width=\linewidth]{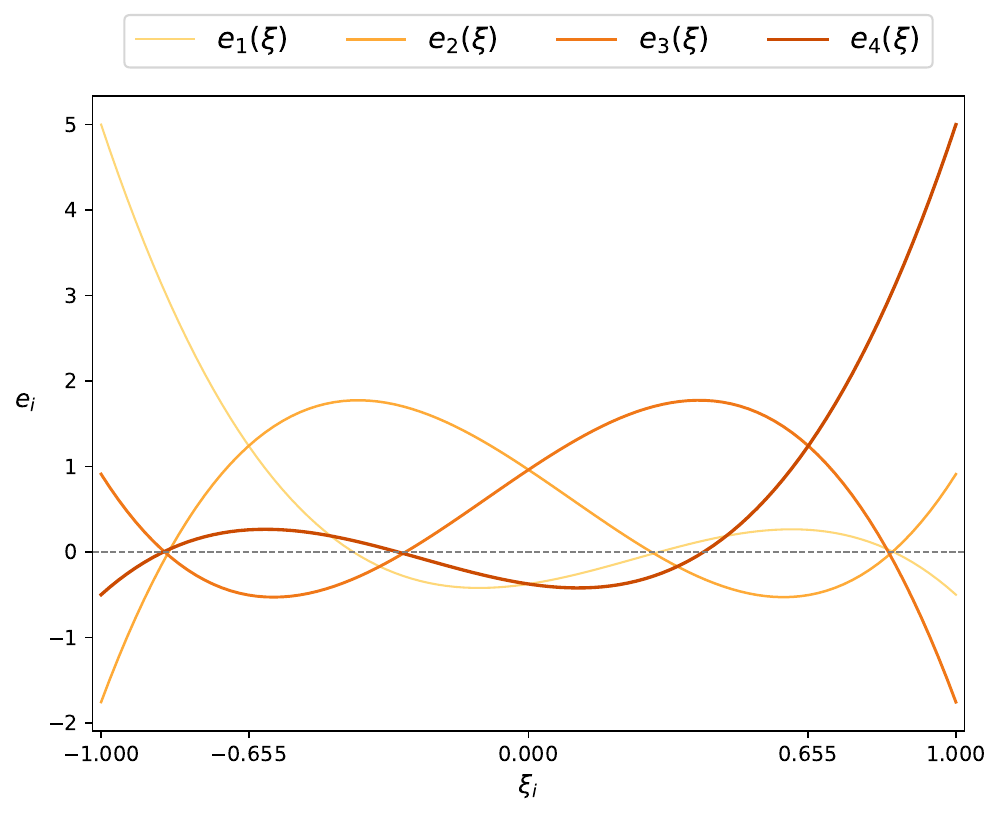}
        \caption{Edge polynomial basis}
    \end{subfigure}
    \caption{ 1D basis functions for polynomial degree N = 4 }
    \label{fig:1d_basis}
\end{figure}

In higher dimensions, basis functions are built by taking tensor products of the one-dimensional nodal and edge polynomials introduced above. 

\subsection{Incidence matrices}
\label{sec:discrete_operators}


The relation between a nodal and an edge expansions becomes clear when we take the derivative of a functions expanded in a nodal basis, \cite{gerritsma2011edge}
\begin{equation}
  \phi_h(\xi) = \sum_{i=0}^N \phi_i h_i(\xi) \quad \Longrightarrow \quad \mathrm{d} \phi_h(\xi) 
  = \sum_{j=1}^N (\phi_j-\phi_{j-1}) e_j(\xi)\; .
  \label{eq:derivative_primal}
\end{equation}
Going from a nodal representation to its derivatives involves two steps: The expansion coefficients of the derivative are the difference of consecutive expansion coefficients of the nodal expansion. And the nodal basis that span the space of polynomials of degree $N$ changes to edge polynomials that span the space of polynomials of degree $(N-1)$. If the expansion coefficients $\phi_i$ are stored in a column vector $\overline{\phi}$, then we can compute the expansion coefficients for the derivative as

\begin{equation}
    \delta\overline{\phi} = \underbrace{\left ( \begin{array}{ccccccc}
    -1 & 1 & 0 & \cdots  & & & \\
    0 & -1 & 1 & & & & \\
    \vdots & & \ddots & \ddots & & & \vdots \\
      & & & & -1 & 1 & 0 \\
      & & & \cdots & 0 & -1 & 1
     \end{array} \right )}_{\mathbb{E}} \overline{\phi} \;.
\end{equation}
Here\footnote{This $\delta$ should not be confused with the $\delta$ from Section~\ref{sec:Largrangina_formulation} where it denotes variations.} $\delta$ denotes the difference operator which acts on the degrees of freedom and is represented by the incidence matrix $\mathbb{E}$. Incidence matrices are coordinate- and metric-free.


For example, consider that a scalar function $\phi$ is discretized on a 2D primal mesh using a nodal ($0$-form) expansion. The approximate function $\phi_h(\xi,\eta)$ is written as
\begin{equation}
  \phi_h(\xi,\eta)
  \;=\;
  \sum_{i=0}^N \sum_{j=0}^N
  \phi_{i,j}\,
  h_i(\xi)\,h_j(\eta),
  \label{eq:phi0}
\end{equation}
where, $h_i(\xi)$ and $h_j(\eta)$ denote the nodal basis functions in the respective directions,
and $\phi_{i,j}$ represent the nodal degrees of freedom. Using \eqref{eq:derivative_primal} we have that
\[ \frac{\mathrm{d}}{\mathrm{d}\xi} \phi_h(\xi,\eta)
  \;=\; \sum_{i=1}^N \sum_{j=0}^N \left ( \phi_{i,j} - \phi_{i-1,j} \right ) e_i(\xi) h_j(\eta) \;, \qquad 
  \frac{\mathrm{d}}{\mathrm{d}\eta} \phi_h(\xi,\eta)
  \;=\; \sum_{i=0}^N \sum_{j=1}^N \left ( \phi_{i,j} - \phi_{i,j-1} \right ) h_i(\xi) e_j(\eta) \;.
\]
Combining these two operations for $N=1$ gives
\begin{equation}
    \nabla \phi_h(\xi,\eta) = \left [ \begin{array}{cccc}
    e_1(\xi)h_0(\eta) & e_1(\xi) h_1(\eta) & 0 &  0 \\
    0 &  0 & h_0(\xi)e_1(\eta)  & h_1(\xi) e_1(\eta)
    \end{array}\right ]
    \underbrace{\left ( \begin{array}{cccc}
    -1 & 1 & 0 & 0 \\
    0 & 0 & -1 & 1 \\
    -1 & 0 & 1 & 0 \\
    0 & -1 & 0 & 1
    \end{array} \right )}_{\mathbb{E}^{1,0}}
    \left ( \begin{array}{c}
    \phi_{1,1} \\
    \phi_{2,1} \\
    \phi_{1,2} \\
    \phi_{2,2}
    \end{array} \right ) \;.
\end{equation}

This matrix, $\mathbb{E}^{1,0}$, is the discrete gradient applied to the degrees of freedom. Similarly, the incidence matrix for the curl operator can be constructed, which gives the incidence matrix $\mathbb{E}^{2,1}$, see \cite{gerritsma2011edge}. Then one can confirm that
\begin{equation*}
  \mathbb{E}^{2,1}\,\mathbb{E}^{1,0} = 0,
  \quad\text{and more generally}\quad
  \mathbb{E}^{k+1,k}\,\mathbb{E}^{k,\,k-1} = 0,
\end{equation*}
which mirror the continuous identities \(\nabla \times (\nabla \varphi) = 0\) and \(\nabla \cdot (\nabla \times \mathbf{u}) = 0\).

\subsection{Dual Basis Functions}


When constructing discrete differential operators, it is advantageous to work with a \emph{dual} basis that is bi-orthogonal to the primal basis where we treat $L^2(\Omega)$ as the pivot space. The dual basis is obtained by a transformation involving the mass matrix. The general formula for obtaining the dual basis function is given by 

\begin{equation}
\tilde{\Psi}^{(n-k)} = \Psi^k\; \left[ \mathbb{M}^k \right]^{-1} ,
\end{equation}

where $\Psi^k$ represents the row vector of basis functions, and $\tilde{\Psi}^{(n-k)}$ denotes its dual. 
Now that we have constructed the dual basis functions, we can also set up the dual mass matrix, which turns out to be
\begin{equation}
    \widetilde{\mathbb{M}}^{n-k} = \int_\Omega \tilde{\Psi}^{(n-k)} \wedge \star \tilde{\Psi}^{(n-k)} = \left[ \mathbb{M}^k \right]^{-1} \left [ \int_\Omega \Psi^k \wedge \star \Psi^k \right ] \left[ \mathbb{M}^k \right]^{-1} = \left[ \mathbb{M}^k \right]^{-1} \;.
    \label{eq:dual_mass_matrix}
\end{equation} 
Note that by construction we have
\begin{equation}
    \int_\Omega \Psi^k \wedge \star \tilde{\Psi}^{(n-k)} = \mathbb{I} \;.
\end{equation}
This implies that duality pairing between a discrete primal $k$-form and dual $(n-k)$-form has a "mass matrix" equal to the identity matrix. This operation does no longer depend explicitly on the basis functions and duality pairing is just the vector product of the expansion coefficients of both forms.
For a more rigorous derivation of dual basis functions and their properties, refer to \cite{jain_gerritsma_quads_hexes}.

\subsubsection{1D Dual Basis Functions}
\label{sec:dual_basis_functions}

In 1D, the dual basis functions are constructed from the primal nodal, $h(\xi)$, and edge, $e(\xi)$, bases using their respective mass matrices.The dual nodal basis $\tilde{h}(\xi)$, defined on the dual grid nodes, is dual to the primal edge functions $e(\xi)$, and conversely the dual edge basis $\tilde{e}(\xi)$, defined on the dual grid edges, is dual to the primal nodal functions $h(\xi)$,

\begin{equation}
  \label{eq:dual_basis_definitions}
  \tilde{h}(\xi) = e(\xi)\; \left[ \mathbb{M}^{(1)} \right]^{-1}\; ,
  \qquad
  \tilde{e}(\xi) = h(\xi)\; \left[ \mathbb{M}^{(0)} \right]^{-1}\; .
\end{equation}

where $\mathbb{M}^{(1)}$ is the primal 1-form (edge) mass matrix and $\mathbb{M}^{(0)}$ is the primal 0-form (nodal) mass matrix.

\begin{equation}
  \label{eq:mass_matrices}
  \mathbb{M}_{ij}^{(1)} = \int\limits_{-1}^{1} e_i(\xi)e_j(\xi)\;\mathrm{d}\xi \quad \text{for } i = 1,\dots,N,
  \qquad
  \mathbb{M}_{ij}^{(0)} = \int\limits_{-1}^{1} h_i(\xi)h_j(\xi)\;\mathrm{d}\xi \quad \text{for } i = 0,\dots,N.
\end{equation}

\begin{figure}[h!]
    \centering
    \begin{subfigure}[b]{0.35\textwidth}
        \centering
        \includegraphics[width=\linewidth]{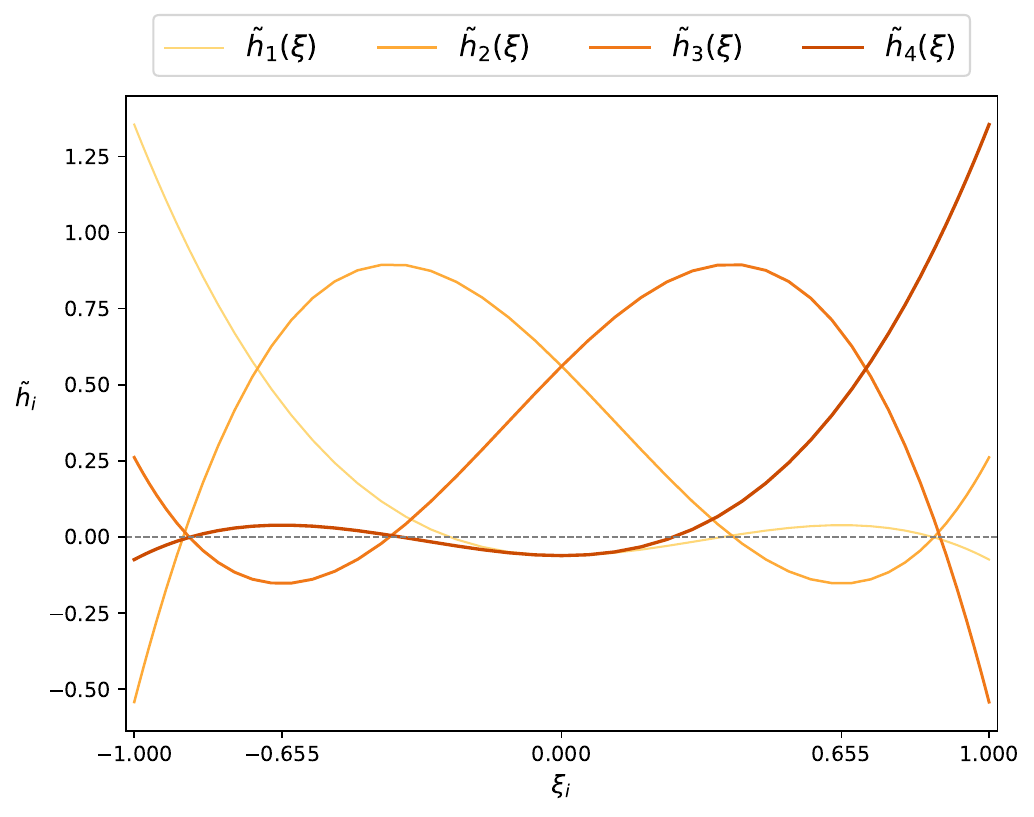}
        \caption{Dual Nodal polynomial basis}
        
    \end{subfigure}
    \hspace{2cm}
    \begin{subfigure}[b]{0.35\textwidth}
        \centering
        \includegraphics[width=\linewidth]{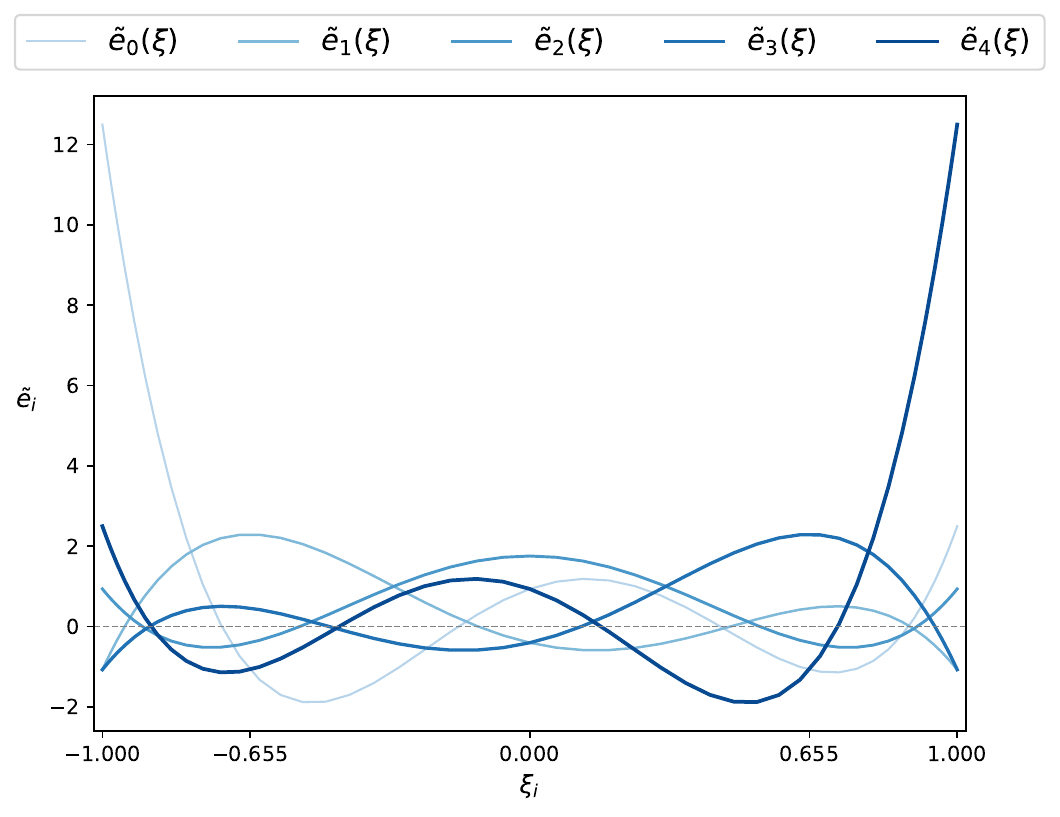}
        \caption{Dual Edge polynomial basis}
        
    \end{subfigure}
    \caption{ 1D Dual basis functions for polynomial degree N = 4 }
    \label{fig:1d_dual_basis}
\end{figure}

\section{Variable Representation in Two Dimensional Space and Time $\mathbb{R}^{2+1}$}\label{sec:SpaceTimeFormulation}

We will apply the mimetic spectral element method to formulation derived in Section~\ref{sec:Largrangina_formulation}. If a two-dimensional Cartesian space is considered, all physical variables are described as \emph{bundle valued differential forms} on a single space-time manifold $\boldsymbol{M} = (\mathcal{B} \times T), \in \mathbb{R}^{2+1}$. This approach allows for a consistent discretization of all components, both spatial and temporal on the manifold's discrete counterpart. As detailed in Table~\ref{tab:variable_representation}, these fields naturally decompose into $x$, $y$ and $t$ components in the spatial and temporal domain. This involves fields that are tensor products.
The fibre space can be scalar-valued, vector-valued (e.g., for flow map $\bm{\varphi}$) or covector-valued (e.g., for momentum $\widetilde{\bm{\pi}}$). In the representations that follow, a clear distinction is maintained between the value space and the domain of the forms. The bold symbols $\bm{\partial}_\alpha$ and $\bm{d}x_\alpha$ denote the basis vectors and covectors, respectively, corresponding to the physical components. In contrast, the regular font differentials $d\xi, d\eta,$ and $d\tau$ represent the fundamental coordinate $1$-forms indicating the directions of integration. Table \ref{tab:variable_representation} provides a detailed summary of this space-time representation for each primary variable. For a more detailed review of the representation of physical variables as differential forms, one can refer to \cite{flanders1989differential,frankel2011geometry}.

{\bf Note on notation}: To distinguish between the kinematic (primal) degrees of freedom and the dynamic (dual) degrees of freedom, we will denote the dual variables with a tilde, $\widetilde{\cdot}$ and the vector of the dual degrees of freedom by $\overline{\widetilde{\cdot}}$, while the vector of primal degrees of freedom will denoted by $\overline{\cdot}$.

\newcolumntype{S}{>{\raggedright\arraybackslash}p{2cm}}
\newcolumntype{C}{>{\centering\arraybackslash}p{1.6cm}}
\newcolumntype{Y}{>{\raggedright\arraybackslash}X} 

\begin{table}[H]
\centering
\small

\renewcommand{\arraystretch}{2.0}
\setlength{\tabcolsep}{4pt}
\begin{tabularx}{\textwidth}{
    S 
    C 
    C 
    C 
    >{\hsize=21.2\hsize\centering\arraybackslash}Y 
    C 
}
\toprule
Variable & Symbol & Bundle type & Degree (space,time) & Representation on $\boldsymbol{M}$ & Placement \\
\midrule
Flow map &
$\bm{\varphi}$ &
vector &
$(0,0)$ &
\adjustbox{max width=14.5\linewidth}{$
\bm{\varphi}^{(0)}=\bm{\partial}_x\!\otimes\!\varphi_x(\xi,\eta,\tau)+\bm{\partial}_y\!\otimes\!\varphi_y(\xi,\eta,\tau)
$} &
primal
\\ 

Velocity &
$\bm{V}$ &
vector &
$(0,1)$ &
\adjustbox{max width=15\linewidth}{$
\bm{V}^{(1)}=\big(\bm{\partial}_x\!\otimes\!V_x(\xi,\eta,\tau)+\bm{\partial}_y\!\otimes\!V_y(\xi,\eta,\tau)\big)\,\mathrm{d}\tau
$} &
primal
\\ 

Deformation gradient &
$\bm{F}$ &
vector &
$(1,0)$ &
\adjustbox{max width=20\linewidth}{$
\bm{F}^{(1)}=\sum_{\alpha\in\{x,y\}}\bm{\partial}_\alpha\!\otimes\!\big(F_{\alpha x}(\xi,\eta,\tau)\ \mathrm{d}\xi+F_{\alpha y}(\xi,\eta,\tau)\ \mathrm{d}\eta\big)
$} &
primal
\\ 

Momentum &
$\widetilde{\bm{\pi}}$ &
covector &
$(2,0)$ &
\adjustbox{max width=19.5\linewidth}{$
\widetilde{\bm{\pi}}^{(2)}=\mathbf{d}x\!\otimes\!\pi_x(\xi,\eta,\tau)\ \mathrm{d}\xi\wedge \mathrm{d}\eta+\mathbf{d}y\!\otimes\!\pi_y(\xi,\eta,\tau)\  \mathrm{d}\xi\wedge \mathrm{d}\eta
$} &
dual
\\ 

First Piola–Kirchhoff stress &
$\widetilde{\bm{P}}$ &
covector &
$(1,1)$ &
\adjustbox{max width=22.5\linewidth}{$
\widetilde{\bm{P}}^{(2)}=\sum_{\alpha\in\{x,y\}}\mathbf{d}x_{\alpha}\!\otimes\!\big(P_{\alpha x}(\xi,\eta,\tau)\ \mathrm{d}\eta-P_{\alpha y}(\xi,\eta,\tau)\ \mathrm{d}\xi\big)\wedge \mathrm{d}\tau
$} &
dual
\\ 
\bottomrule
\end{tabularx} 
\vspace{0.2cm}
\caption{Bundle valued space–time form representation of physical variables.}
\label{tab:variable_representation}
\end{table}

\subsection{Primal-Dual Variable Staggering}

\subsubsection{Configuration (Kinematic) Variables:}
\begin{figure}[h!]
\centering
  \includegraphics[width=10.5cm]{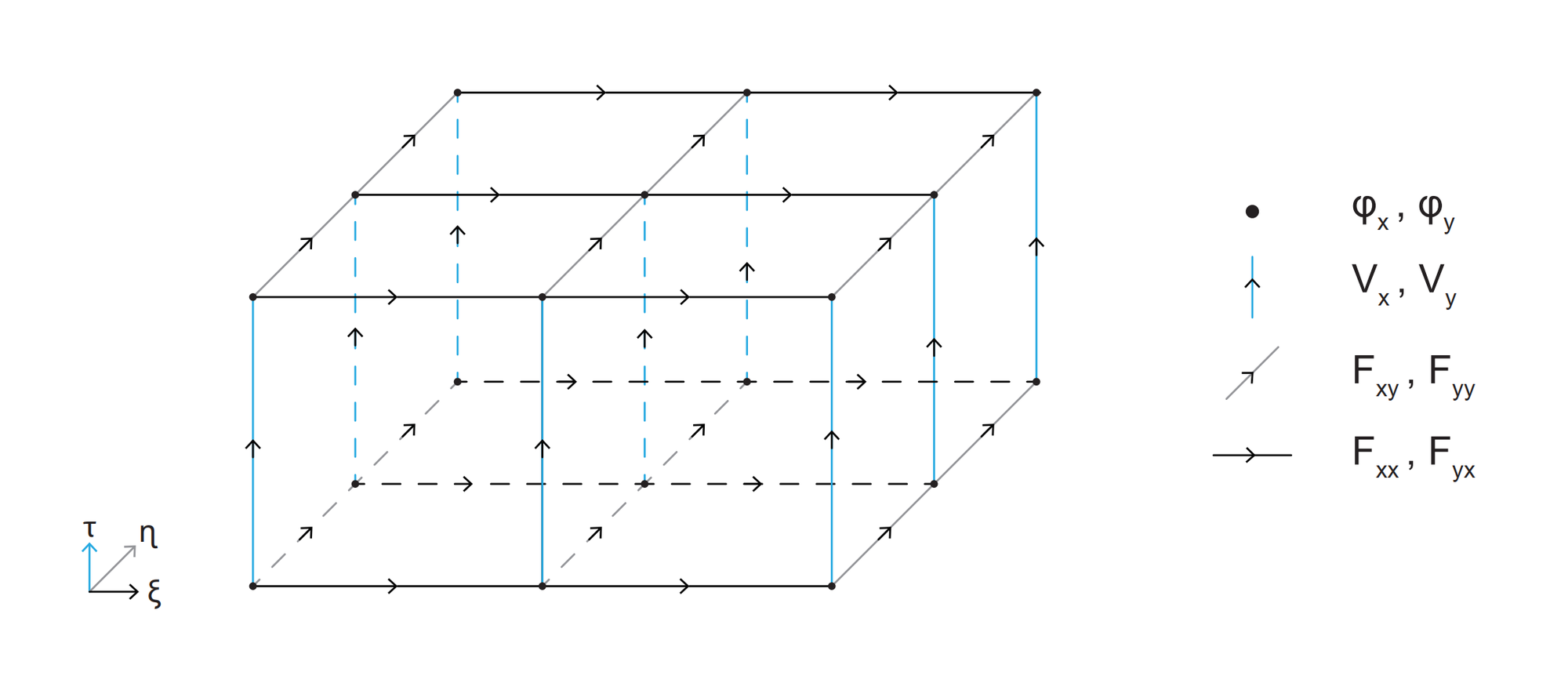}
  \caption{Localization of degrees of freedom on the primal mesh.}
  \label{fig:primal}
\end{figure}

In a Lagrangian framework, the fields defining the system's geometry are the kinematic variables, the flow map $\boldsymbol{\varphi}$, velocity $\boldsymbol{V} $ and deformation gradient $\boldsymbol{F}$, and they appear naturally in the Lagrangian density \eqref{eq:lagrangian}. These vector valued fields are associated with the primal mesh, the direct discrete representation of the domain. As shown in Figure \ref{fig:primal} the flow map components ($\varphi_x, \varphi_y$) are localized at the primal nodes (vertices). The degrees of freedom for velocity and deformation gradient are defined as integrals the edges, which are equivalent to the difference of the flow map degrees of freedom at the edge's endpoints, the differences in the temporal directions are assigned to the vertical edges in the space-time mesh, as indicated in the figure and the differences in the spatial plane are assigned to the horizontal edges.
    
\begin{equation}\label{eq:edge-integrals}
\begin{aligned}
\bar{V}_{x,i,j} =\varphi_{x,i,j}^{n+1}-\varphi_{x,i,j}^{n} &= \int\limits_{t^n}^{t^{n+1}} V_x(x_i,y_j,s)\,ds,
& \text{and}
\qquad 
\bar{V}_{y,i,j} =\varphi_{y,i,j}^{n+1}-\varphi_{y,i,j}^{n} &= \int\limits_{t^n}^{t^{n+1}} V_y(x_i,y_j,s)\,ds,\\[0.3em]
\bar{F}_{xx,i,j}= \varphi_{x,i,j}^{n}-\varphi_{x,i-1,j}^{n} &= \int\limits_{x_{i-1}}^{x_i} F_{xx}(s,y_j,t^{n})\,ds,
&\text{and}\qquad \bar{F}_{yx,i,j}=
\varphi_{y,i,j}^{n}-\varphi_{y,i-1,j}^{n} &= \int\limits_{x_{i-1}}^{x_i} F_{yx}(s,y_j,t^{n})\,ds,\\[0.3em]
\bar{F}_{xy,i,j}=\varphi_{x,i,j}^{n}-\varphi_{x,i,j-1}^{n} &= \int\limits_{y_{j-1}}^{y_j} F_{xy}(x_i,s,t^{n})\,ds,
&\text{and}\qquad \bar{F}_{yy,i,j }=
\varphi_{y,i,j}^{n}-\varphi_{y,i,j-1}^{n} &= \int\limits_{y_{j-1}}^{y_j} F_{yy}(x_i,s,t^{n})\,ds.
\end{aligned}
\end{equation}

\subsubsection{Source (Dynamic) Variables:}  

\begin{figure}[h!]
    \centering
    \begin{minipage}{0.35\textwidth}
        \centering
        \includegraphics[width=\linewidth]{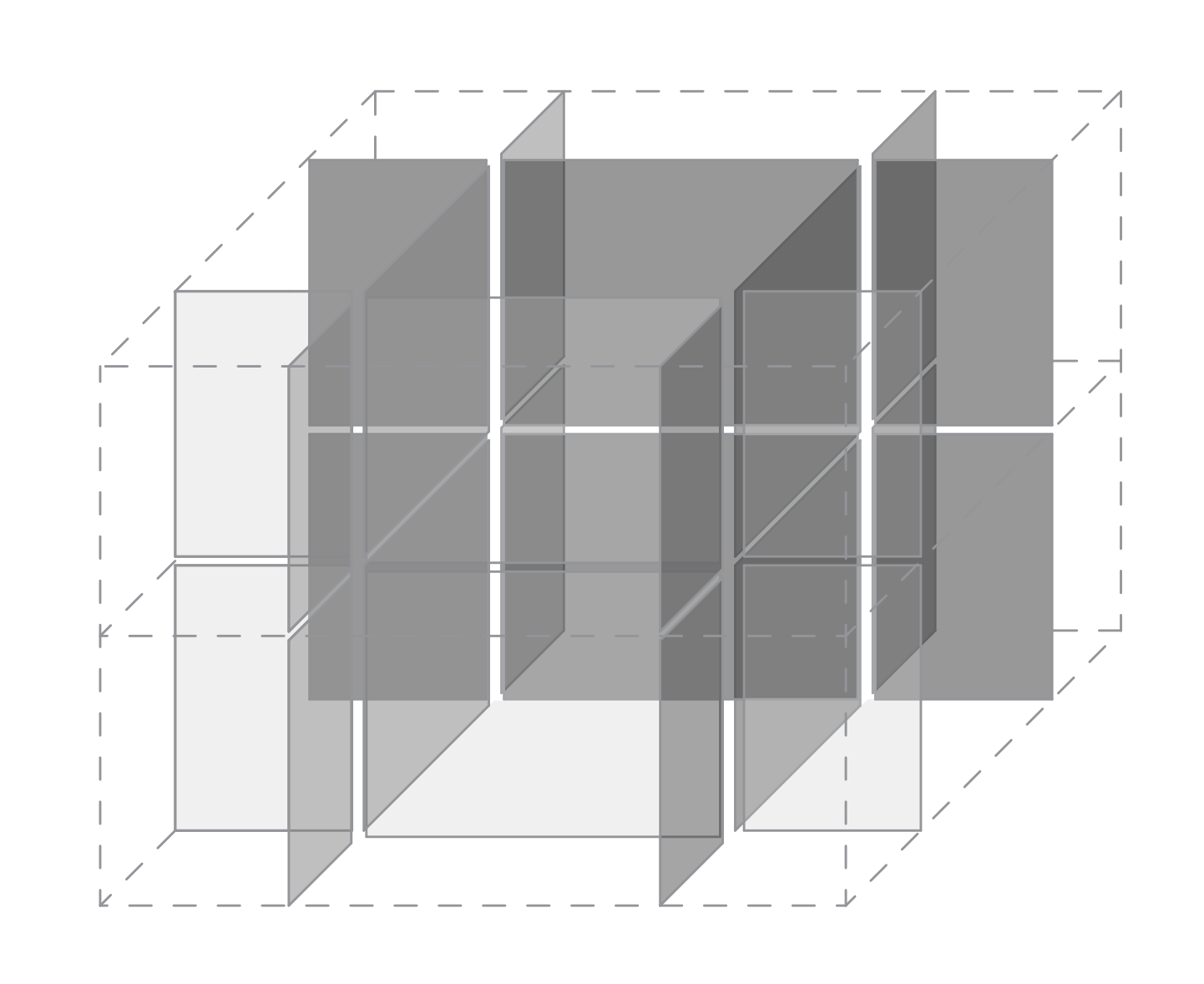}
        \subcaption{First Piola-Kirchhoff stress degrees of freedom on dual spatio-temporal faces.}
    \end{minipage}
    \hfill
    \begin{minipage}{0.27\textwidth}
        \centering
        \includegraphics[width=\linewidth]{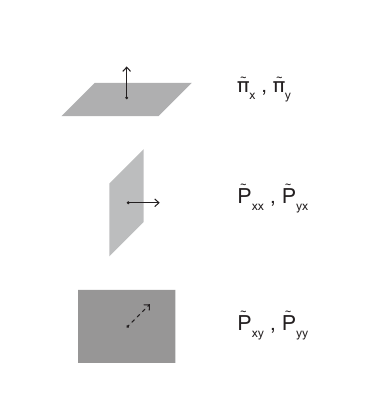}
    \end{minipage}
    \hfill
    \begin{minipage}{0.32\textwidth}
        \centering
        \includegraphics[width=\linewidth]{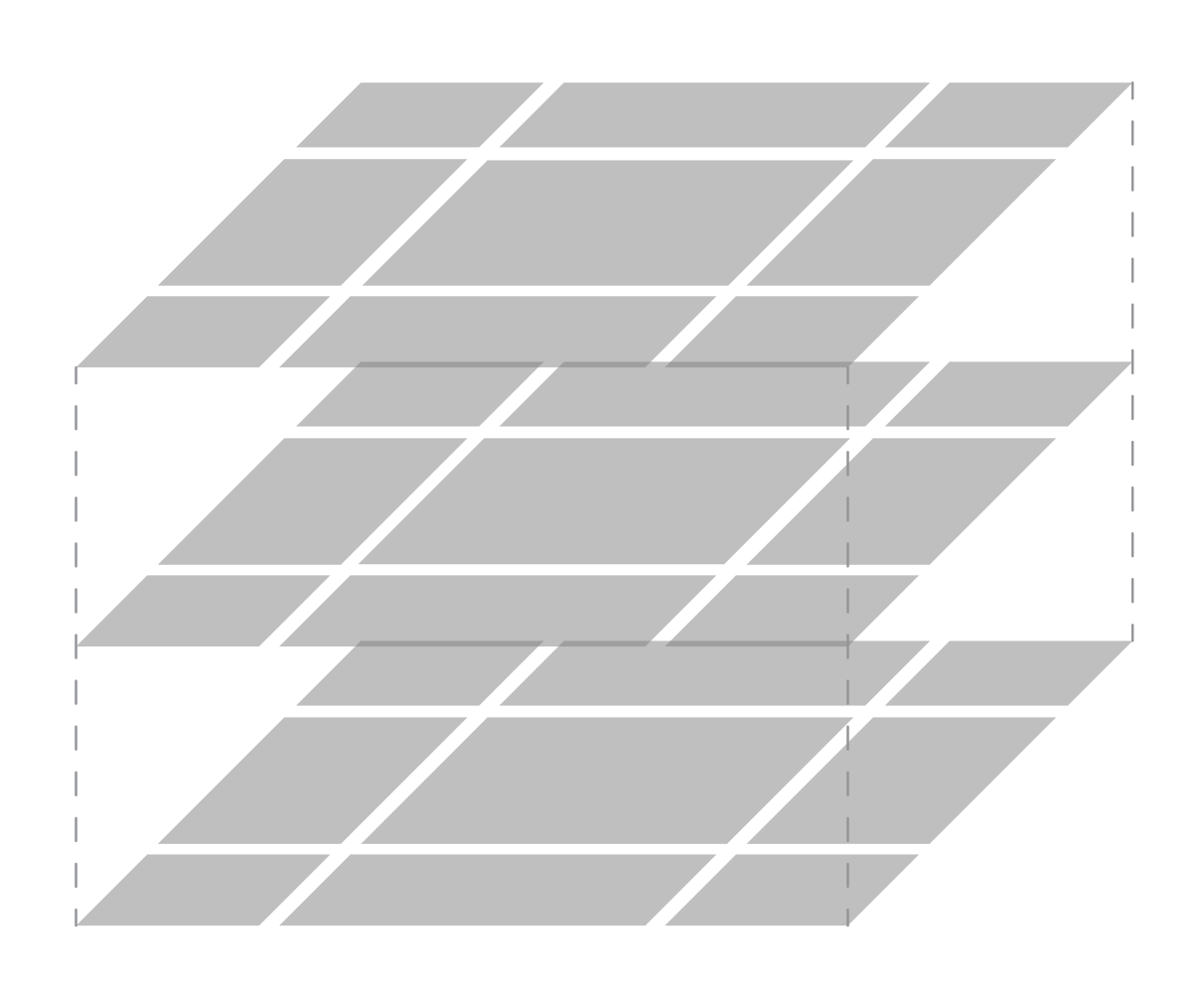}
        \subcaption{Momentum degrees of freedom on dual spatial faces.}
    \end{minipage}
    \caption{Localization of the degrees of freedom on the dual mesh.}
    \label{fig:dual_grid}
\end{figure}

By contrast, covector valued variables (forces, momentum, stresses) act on the configuration and could be thought of as being placed on the dual grid as shown in Figure \ref{fig:dual_grid}. Although the dual cells and source degrees of freedom are not localized as on the primal grid, their action is visible through the evolution of the primal configuration. Here momentum is represented on dual spatial surfaces, these faces are naturally paired with the temporal edges of the primal mesh that carry the velocity. The first Piola–Kirchhoff stress degrees of freedom are associated with dual spatio-temporal faces. These stress degrees of freedom represent the impulse of the surface forces (i.e., stress integrated over a spatial edge and a dual time interval). These dual faces are naturally paired with the primal edges that carry the deformation gradient, to produce kinetic energy or the work done by the surface forces.

\subsubsection{Spectral Basis Expansion}

For the \(\mathbb{R}^{(2+1)}\) grid, basis functions are constructed at the Gauss-Lobatto-Legendre (GLL) co-ordinates $\xi \in \{\xi_0, \;\xi_1, \dots, \xi_{N}\}, \quad  
\eta \in \{\eta_0, \;\eta_1, \dots, \eta_{N}\}, \quad  
\tau \in \{\tau_0, \;\tau_1, \dots, \tau_{N_t}\}$ via the tensor product of the 1D basis functions presented earlier, as shown below.

\paragraph{Primal Nodal Basis:}
The primal nodal basis functions \(\boldsymbol{\Psi}_0\), are defined as the tensor product of one-dimensional nodal Lagrange basis functions, and belong to the discrete space,  which can also be interpreted as the space \( H^1(\mathcal{B} \times T) \subset L^2(\mathcal{B} \times T)\).

\begin{equation}
    \boldsymbol{\Psi}_0 = h(\xi) \otimes h(\eta) \otimes h(\tau) 
\end{equation}

which, in index notation, expands to \(\boldsymbol{\Psi}_0 = h_i(\xi) \;h_j(\eta)\; h_k(\tau), \quad \forall \; i,j = 0,1,\dots,N, \quad k = 0,1,\dots,N_t.\) \\

\paragraph{Primal Edge Basis:}

The edge basis are obtained similarly by tensor product of combination of edge basis in that respective direction and nodal basis in remaining directions. They belong to the discrete space \(\Lambda_h^{(1)}(\mathcal{B} \times T)\) which can also be interpreted as the space \(H\;(\mathrm{curl}\;;\;\mathcal{B} \times T)\). In fact, application of the curl to this bundle-valued form produced the Piola identity, which just follows from curl$\cdot$grad$\equiv 0$.\\

\begin{equation}
    \boldsymbol{\Psi}_1 =\begin{bmatrix}
e(\xi)\otimes h(\eta)\otimes h(\tau) & 0 & 0 \\
0 & h(\xi)\otimes e(\eta) \otimes h(\tau) & 0 \\
0 & 0 & h(\xi)\otimes h(\eta)\otimes e(\tau) 
\end{bmatrix} \\
\end{equation}

\vspace{0.2cm}

\begin{itemize}
    \item Edge in the \(\xi\)-direction : \(\boldsymbol{\Psi}_1^{\xi} = e_i(\xi) \;h_j(\eta)\; h_k(\tau), \quad \forall \; i=1,2,\dots,N, \quad j = 0,1,\dots,N, \quad k = 0,1,\dots,N_t.\)
    \item Edge in the \(\eta\)-direction : \(\boldsymbol{\Psi}_1^{\eta} = h_i(\xi) \;e_j(\eta)\; h_k(\tau), \quad \forall \; i=0,1,\dots,N, \quad j = 1,2,\dots,N, \quad k = 0,1,\dots,N_t.\)
    \item Edge in the \(\tau\)-direction : \(\boldsymbol{\Psi}_1^{\tau} = h_i(\xi) \;h_j(\eta)\; e_k(\tau), \quad \forall \; i,j=0,1,\dots,N,  \quad k = 1,2,\dots,N_t.\)
\end{itemize}


\paragraph{Dual Surface Basis:}

The dual surface basis functions are constructed by taking the tensor product of a dual edge basis function in two coordinate directions and a nodal basis in the perpendicular direction. They belong to the discrete space \(\tilde\Lambda_h^{(2)}(\mathcal{B} \times T)\) which can also be interpreted as the space \(H\;(div\;;\;\mathcal{B} \times T)\).

\begin{equation}
    \boldsymbol{\tilde{\Psi}}_2 =
    \begin{bmatrix}
        \tilde{h}(\xi) \otimes \tilde{e}(\eta) \otimes \tilde{e}(\tau) & 0 & 0 \\
        0 & \tilde{e}(\xi) \otimes \tilde{h}(\eta) \otimes \tilde{e}(\tau) & 0 \\
        0 & 0 & \tilde{e}(\xi) \otimes \tilde{e}(\eta) \otimes \tilde{h}(\tau)
    \end{bmatrix}
\end{equation}

\begin{itemize}
    \item Surface \(\perp\) to the \(\xi\)-direction :  $\boldsymbol{\Psi}_2^{\xi} = \tilde{h}_i(\xi) \, \tilde{e}_j(\eta) \, \tilde{e}_k(\tau), 
    \quad \forall \; i=0,1,\dots,N, \quad j = 1,2,\dots,N, \quad \\ k = 1,2,\dots,N_t.$
    
    \item Surface \(\perp\) to the \(\eta\)-direction : $    \boldsymbol{\Psi}_2^{\eta} = \tilde{e}_i(\xi) \, \tilde{h}_j(\eta) \, \tilde{e}_k(\tau), 
    \quad \forall \; i=1,2,\dots,N, \quad j = 0,1,\dots,N, \quad \\ k = 1,2,\dots,N_t.$
   
    \item Surface \(\perp\) to the \(\tau\)-direction : $\boldsymbol{\Psi}_2^{\tau} = \tilde{e}_i(\xi) \, \tilde{e}_j(\eta) \, \tilde{h}_k(\tau), 
    \quad \forall \; i=1,2,\dots,N, \quad j = 1,2,\dots,N, \quad \\ k =0,1,\dots,N_t.$
    
\end{itemize}

Based on these definitions, the discrete physical variables are expanded as shown in Tables~\ref{tab:spec_exp_kin} and \ref{tab:spec_exp_dyn}.

\newcolumntype{C}[1]{>{\centering\arraybackslash}m{#1}}

\begin{table}[H]
\centering
\renewcommand{\arraystretch}{1.4} 
\begingroup

\setlength{\abovedisplayskip}{4pt}
\setlength{\belowdisplayskip}{4pt}
\setlength{\abovedisplayshortskip}{2pt}
\setlength{\belowdisplayshortskip}{2pt}
\begin{tabular}{|C{2.0cm}|C{2.0cm}|C{1.5cm}|C{1.6cm}|C{6.0cm}|}
\hline
\textbf{Variable} & \textbf{Association} & \textbf{Basis} & \textbf{DOF's} & \textbf{Expansion} \\ \hline

\textbf{Flowmap} $\boldsymbol{\varphi}$ 
& Primal nodes 
& $\boldsymbol{\Psi}_0$ 
& $\overline{\boldsymbol{\varphi}}_x$, $\overline{\boldsymbol{\varphi}}_y$ 
& {\small
\[
\begin{aligned}
\varphi_x^h(\xi,\eta,\tau) &= \sum_{i,j,k=0}^{N,N,N_t} 
{\varphi}_{ijk}^x\, h_i(\xi)\, h_j(\eta)\, h_k(\tau) \\
\varphi_y^h(\xi,\eta,\tau) &= \sum_{i,j,k=0}^{N,N,N_t} 
{\varphi}_{ijk}^y\, h_i(\xi)\, h_j(\eta)\, h_k(\tau)
\end{aligned}
\]
} \\ \hline

\textbf{Deformation Gradient} $\bm{F}$ 
& Primal spatial edges 
& $\boldsymbol{\Psi}_1^{\xi}$, $\boldsymbol{\Psi}_1^{\eta}$ 
& $\overline{\bm{F}}_{xx}$, $\overline{\bm{F}}_{xy}$ $\overline{\bm{F}}_{yx}$, $\overline{\bm{F}}_{yy}$ 
& {\small
\[
\begin{aligned}
F_{xx}^h &= \sum_{i=1}^{N}\sum_{j,k=0}^{N,N_t} 
{F}_{ijk}^{xx}\, e_i(\xi)\, h_j(\eta), h_k(\tau) \\
F_{xy}^h &= \sum_{i=0}^{N}\sum_{j=1}^{N}\sum_{k=0}^{N_t} 
{F}_{ijk}^{xy}\, h_i(\xi)\, e_j(\eta)\, h_k(\tau) \\
F_{yx}^h &= \sum_{i=1}^{N}\sum_{j=0}^{N}\sum_{k=0}^{N_t} 
{F}_{ijk}^{yx}\, e_i(\xi)\, h_j(\eta)\, h_k(\tau) \\
F_{yy}^h &= \sum_{i=0}^{N}\sum_{j=1}^{N}\sum_{k=0}^{N_t} 
{F}_{ijk}^{yy}\, h_i(\xi)\, e_j(\eta)\, h_k(\tau)
\end{aligned}
\]
} \\ \hline

\textbf{Velocity} $\bm{V}$ 
& Primal (temporal) edges 
& $\boldsymbol{\Psi}_1^{\tau}$ 
& $\overline{\bm{V}}_x$, $\overline{\bm{V}}_y$ 
& {\small
\[
\begin{aligned}
V_x^h &= \sum_{i,j,k=0}^{N,N,N_t}
{V}_{ijk}^x\, h_i(\xi)\, h_j(\eta)\, h_k(\tau) \\
V_y^h &= \sum_{i,j=0}^{N,N}\sum_{k=1}^{N_t}
{V}_{ijk}^y\, h_i(\xi)\, h_j(\eta)\, e_k(\tau)
\end{aligned}
\]
} \\ \hline

\end{tabular}
\endgroup 
\vspace{0.5cm}
\caption{Summary of field associations, basis functions, and expansions for $\boldsymbol{\varphi}$, $\bm{F}$, and $\bm{V}$.}
\label{tab:spec_exp_kin}
\end{table}

\begin{table}[H]
\centering
\renewcommand{\arraystretch}{1.4} 
\begingroup

\setlength{\abovedisplayskip}{4pt}
\setlength{\belowdisplayskip}{4pt}
\setlength{\abovedisplayshortskip}{2pt}
\setlength{\belowdisplayshortskip}{2pt}
\begin{tabular}{|C{2.0cm}|C{2.0cm}|C{1.5cm}|C{1.6cm}|C{6.0cm}|}
\hline
\textbf{Variable} & \textbf{Association} & \textbf{Basis} & \textbf{DOF's} & \textbf{Expansion} \\ \hline

\textbf{Momentum} $\widetilde{\bm{\pi}}$ 
& Dual surfaces 
& $\widetilde{\boldsymbol{\Psi}}_2^{\tau}$ 
& $\overline{\widetilde{\bm{\pi}}}_x$, $\overline{\widetilde{\bm{\pi}}}_y$ 
& {\small
\[
\begin{aligned}
\widetilde{\pi}_x^h &= \sum_{i,j=1}^{N,N}\sum_{k=0}^{N_t}
\widetilde{\pi}^x_{ijk}\, \widetilde{e}_i(\xi)\, \widetilde{e}_j(\eta)\, \widetilde{h}_k(\tau) \\
\widetilde{\pi}_y^h &= \sum_{i,j=1}^{N,N}\sum_{k=0}^{N_t}
\widetilde{\pi}^y_{ijk}\, \widetilde{e}_i(\xi)\, \widetilde{e}_j(\eta)\, \widetilde{h}_k(\tau)
\end{aligned}
\]
} \\ \hline

\textbf{First Piola--Kirchhoff Stress} $\widetilde{\bm{P}}$ 
& Dual spatio-temporal surfaces 
& $\widetilde{\boldsymbol{\Psi}}_2^{\xi}$, $\widetilde{\boldsymbol{\Psi}}_2^{\eta}$ 
& $\overline{\widetilde{\bm{P}}}_{xx}$, $\overline{\widetilde{\bm{P}}}_{xy}$ $\overline{\widetilde{\bm{P}}}_{yx}$, $\overline{\widetilde{\bm{P}}}_{yy}$ 
& {\small
\[
\begin{aligned}
\widetilde{P}_{xx}^h &= \sum_{i=0}^{N}\sum_{j,k=1}^{N,N_t}
\widetilde{P}_{ijk}^{xx}\, \widetilde{h}_i(\xi)\, \widetilde{e}_j(\eta)\, \widetilde{e}_k(\tau) \\
\widetilde{P}_{xy}^h &= \sum_{i=1}^{N}\sum_{j=0}^{N}\sum_{k=1}^{N_t}
\widetilde{P}_{ijk}^{xy}\, \widetilde{e}_i(\xi)\, \widetilde{h}_j(\eta)\, \widetilde{e}_k(\tau) \\
\widetilde{P}_{yx}^h &= \sum_{i=0}^{N}\sum_{j=1}^{N}\sum_{k=1}^{N_t}
\widetilde{P}_{ijk}^{yx}\, \widetilde{h}_i(\xi)\, \widetilde{e}_j(\eta)\, \widetilde{e}_k(\tau) \\
\widetilde{P}_{yy}^h &= \sum_{i=1}^{N}\sum_{j=0}^{N}\sum_{k=1}^{N_t}
\widetilde{P}_{ijk}^{yy}\, \widetilde{e}_i(\xi)\, \widetilde{h}_j(\eta)\, \widetilde{e}_k(\tau)
\end{aligned}
\]
} \\ \hline

\end{tabular}
\endgroup 
\vspace{0.5cm}
\caption{Summary of field associations, basis functions, and expansions for $\tilde{\bm{\pi}}$ and $\bm{P}$.}
\label{tab:spec_exp_dyn}
\end{table}

\section{Algebraic Reduction}
\label{sec:algebraic_reduction}

This section details the algebraic reduction of the continuous weak forms to a discrete system of matrix equations. The process begins by replacing continuous fields with their discrete spectral expansions, by replacing the variations with test functions, and replacing continuous operators and inner products with their discrete mimetic counterparts (incidence matrices \(\mathbb{E}\) and Hodge matrices \(\mathbb{H}\)).

\subsection{Discretization of the Momentum Balance Equation}

We begin with the weak form of the momentum balance,~\eqref{eq:final_weak_formulation}, tested with an arbitrary test function. The test function $\phi$ in the discrete setting is chosen  in $H^1(\mathcal{B}\times T)$ thus serves the same purpose as a discrete variation . Note that the tilde symbol over a function indicates that we use a dual representation

\begin{equation}
    \int\limits_{\mathcal{B}\times T } \widetilde{\bm\pi}\;\frac{\partial\phi }{\partial t} \, d(\mathcal{B}\times T) - \int\limits_{\mathcal{B}} \left .\left \langle \widetilde{\bm{\pi}}, \phi \right \rangle \right |_{t=T}\ \mathrm{d} \mathcal{B} + \int\limits_{\mathcal{B}} \left .\left \langle \widetilde{\bm{\pi}}, \phi \right \rangle \right |_{t=0}\ \mathrm{d} \mathcal{B}
    + \int\limits_{\mathcal{B}\times T } \left( (P_{W}-P_{ext} )J\bm{F}^{-T} \right )\; :\nabla\phi\; \, d(\mathcal{B}\times T) = 0\;,
    \label{eq:weak_form_start} 
\end{equation}
for all $\phi \in H^1(\mathcal{B}\times T)$. Here we made the assumption that the external pressure field is constant, which allows us to contract it with the thermodynamic pressure, and we effectively used the gauge pressure with respect to the constant external pressure.
Note that we could have written $\int\limits_{\mathcal{B}\times T } (P_{W}J\bm{F}^{-T})\; :\nabla\phi = \langle \widetilde{\bf{P}}, \nabla \phi \rangle_{\mathcal{B}\times T}$, where $\widetilde{\bf{P}}$ is the first Piola-Kirchhoff stress tensor, but since we do not introduce this variable, we will stick to the formulation in (\ref{eq:weak_form_start}).




\subsubsection{Momentum Time Derivative and Boundary Terms :}

The test function \(\phi\) is expanded in the primal nodes (\(\boldsymbol{\Psi}_0\)) (see Table~\ref{tab:spec_exp_kin}) with degrees of freedom \(\overline{\phi}\). The momentum \(\widetilde{\bm{\pi}}\) is expanded on dual spatial surfaces at a dual time instant (\(\tilde{\boldsymbol{\Psi}}_2^\tau\)) with degrees of freedom \(\overline{\widetilde{\bm{\pi}}}\). Duality pairing \(\left \langle \widetilde{\bm\pi} \;, \frac{\partial \phi}{\partial t} \right \rangle \) for these discretized fields reduces to the vector product of the degrees of freedom, $\overline{\widetilde{\bm{\pi}}}^\top \mathbb{E}_t^{1,0}\overline{\phi} = \left( {\mathbb{E}_t^{1,0}}^\top \overline{\widetilde{\bm{\pi}}} \right )^\top \overline{\phi}$. 
The term 
\[ \int\limits_{\mathcal{B}} \left .\left \langle \widetilde{\bm{\pi}}, \phi \right \rangle \right |^{t=T}_{t=0}\ \mathrm{d} \mathcal{B} = {\overline{\phi}}^T \mathbb{N}_{\pi}^{T} \widehat{\widetilde{\bm{\pi}}}^t - {\overline{\phi}}^T \mathbb{N}_{\pi}^{0} \widehat{\widetilde{\bm{\pi}}}^{t_0} \;,\]
where $\mathbb{N}_{\pi}^{t}$ represents the topological inclusion operator acting on the dual trace variables for momentum. See \cite{jain_gerritsma_quads_hexes} for a discussion on these topological matrices.


\subsubsection{Stress Term :}

This is the non-linear term. The test gradient \(\nabla\phi\) is discretized on primal spatial edges using the spatial gradient operators \(\mathbb{E}_x^{1,0} \ \text{and} \ \mathbb{E}_y^{1,0}\). The First Piola-Kirchhoff stress, \(P_WJ\bm{F}^{-T}\), is related to the deformation gradient \(\bm{F}\) which is discretized via \(\overline{\bm{F}} = \mathbb{E}_{x,y}^{1,0} \ \overline{\bm{\varphi}}\) through the non-linear pressure-weighted mass matrix \(\mathbb{M}_{P_w}\) as shown below.

\begin{equation*}
     \left( P_WJ\bm{F}_x^{-T} \bm: {\nabla\bm{\phi}_x}  \right)_{\mathcal{B} \times T} \;  =\; \overline{\bm\phi}_x^T \;\; \mathbb{E}_x^{1,0^T} \; \mathbb{M}_{P_w}\;\mathbb{E}^{1,0}_y \ \overline{\bm\varphi}_y \; - \;  \overline{\bm\phi}_x^T \;\; \mathbb{E}_y^{1,0^T} \;  \mathbb{M}_{P_w}\;   \mathbb{E}_x^{1,0} \ \overline{\bm\varphi}_y.
\end{equation*}

\begin{equation}
    \left( P_WJ\bm{F}_y^{-T} \bm: {\nabla\bm{\phi}_y}  \right)_{\mathcal{B} \times T} \;  =\; -\; \overline{\bm\phi}_y^T \;\; \mathbb{E}_x^{1,0^T} \; \mathbb{M}_{P_w}\;\mathbb{E}^{1,0}_y \ \overline{\bm\varphi}_x \; + \;  \overline{\bm\phi}_y^T \;\; \mathbb{E}_y^{1,0^T} \;  \mathbb{M}_{P_w}\;   \mathbb{E}_x^{1,0} \ \overline{\bm\varphi}_x.
    \label{eq:discrete_pk}
\end{equation}

In the above equations, discrete inner product \(\overline{{\nabla \bm{\phi}}}^T \mathbb{M}_{P_w} \bar{\bm{F}}\), and thus $ \mathbb{M}_{P_w}$,  is obtained by integrating the product of primal edge basis functions as shown below,

\begin{equation*}
   \mathbb{M}_{P_w} =  \int\limits_{\mathcal{B} \times T} P_w(\rho_0, J) \bigl(\Psi_1^{\xi}\bigr)^{\!T} \,\Psi_1^{\eta} \,d(\mathcal{B} \times T) \;.
\end{equation*}

The term \(P_w(\rho_0, J)\) serves as a spatial and temporal weight, where \(\rho_0\) is the reference density and \(J = \det(\mathbf{F})\) depends on the solution \(\varphi\) by \(\mathbf{F} = \nabla \bm{\varphi}\). Since \(P_w(\rho_0, J)\) depends on \(J\), and thus on \(\varphi\), the matrix \(\mathbb{M}_{P_w}\) becomes highly non-linear. This nonlinearity is handled through a Picard (fixed-point) iteration. Combining these discrete terms derived in the last two sections and factoring out the arbitrary \(\overline{\bm{\phi}}\) gives the final discrete equations upon moving the known initial condition to the right.

\begin{equation}
    \begin{aligned}
         -\mathbb{E}^{{1,0}^T} \overline{\widetilde{\bm\pi}}_x +\;  \mathbb{N}^t_{\pi}\widehat{\widetilde{\bm\pi}}_x^t \; + \; \mathbb{E}_x^{1,0^T} \; \mathbb{M}_{P_w}\;\mathbb{E}^{1,0}_y \ \overline{\bm\varphi}_y \; - \;   \mathbb{E}_y^{1,0^T} \;  \mathbb{M}_{P_w}\;   \mathbb{E}_x^{1,0} \ \overline{\bm\varphi}_y &= - \mathbb{N}_{\pi}^{t_0} \widehat{\widetilde{\bm\pi}}^{t_0}_x  \;, \\ \\
         -\mathbb{E}^{{1,0}^T} \overline{\widetilde{\bm\pi}}_y +\;  \mathbb{N}^t_{\pi}\widehat{\widetilde{\bm\pi}}_y^t \; - \;  \mathbb{E}_x^{1,0^T} \; \mathbb{M}_{P_w}\;\mathbb{E}^{1,0}_y \ \overline{\bm\varphi}_x \; + \;   \mathbb{E}_y^{1,0^T} \;  \mathbb{M}_{P_w}\;   \mathbb{E}_x^{1,0} \ \overline{\bm\varphi}_x &=  - \mathbb{N}_{\pi}^{t_0} \widehat{\widetilde{\bm\pi}}^{t_0}_y \;.
    \end{aligned}
    \label{eq:final_mom_eqn}
\end{equation}

\subsubsection{Discretization of the Constitutive Equation}

The relation between velocity and momentum is given by (\ref{eq:final_weak_formulation}). The test function \(\widetilde{\bm{\lambda}}\) is expanded on the dual spatial surfaces at a constant dual time instant expanded in (\(\widetilde{\boldsymbol{\Psi}}^\tau_2\)) with cochains \(\overline{\widetilde{\bm{\lambda}}}\). The momentum \(\widetilde{\bm{\pi}}\) is expanded on dual spatial surfaces (\(\tilde{\boldsymbol{\Psi}}_2^\tau\)) with cochains \(\overline{\widetilde{\bm{\pi}}}\). 

\begin{equation}
    \int\limits_{\mathcal{B}\times T } \frac{1}{\rho_0}\widetilde{\bm \pi} \;\widetilde{\bm{\lambda}} \ \ \ d(\mathcal{B}\times T) \;
    - \int\limits_{\mathcal{B}\times T } \frac{\partial \varphi}{\partial t}\;\widetilde{\bm{\lambda}} \; \ \ \ d(\mathcal{B}\times T) \; = \; 0 , 
    \quad \quad  \forall \; \widetilde{\bm{\lambda}} \; .
    \label{eq:weak_form}
\end{equation}
The first integral represents an inner product, while the second integral is the duality pairing between primal velocity and dual momentum degrees of freedom.

\begin{equation}
    \widetilde{\bm{\lambda}}^T \left[ \  \int\limits_{\mathcal{B}\times T } \frac{1}{\rho_0} (\widetilde{\Psi}^{\tau}_2)^T\; \widetilde{\Psi}_2^{\tau} \;\; d(\mathcal{B}\times T) \ \right ]\bar{\bm \pi} \; -\;  \widetilde{\bm{\lambda}}^T  \; \mathbb{E}_t^{(1,0)} \;  \bar{\bm \varphi} \; = \; 0.
\end{equation}

The integral within the square brackets is the mass density-weighted mass matrix. Using \eqref{eq:dual_mass_matrix}, this mass matrix is the inverse of the primal mass matrix composed of basis functions for the velocity field \(\mathbb{M}_{\rho_o}\), which maps primal temporal edge degrees of freedom of the velocity to the dual spatial surface degrees of freedom of the momentum, thus the full discrete form of the momentum constitutive equation can be obtained as,

\begin{equation}
    \begin{aligned}
         \mathbb{M}_{\rho_o}^{-1} \overline{\widetilde{\bm \pi}}_x \; - \; \mathbb{E}_t^{1,0} \overline{\bm\varphi}_x \; = \;  0 \\ \\
         \mathbb{M}_{\rho_o}^{-1} \overline{\widetilde{\bm \pi}}_y \; - \; \mathbb{E}_t^{1,0} \overline{\bm\varphi}_y \; = \;  0 
    \end{aligned}
    \label{eq:final_mom_cons_eqn}
\end{equation}

\subsubsection{Initial Conditions and Final System}
\label{sec:initial_conditions_final_system}

To formulate a well-posed problem, we still need initial conditions for the flow map. These two fields are handled using distinct approaches. 
The initial flow map, \(\varphi(\mathbf{X}, t_0) = \widehat{\varphi}^{t_0} \), is prescribed in a strong sense. This is achieved by introducing an additional constraint equation

\begin{equation}
  \int\limits_{\mathcal{B}\times t_0} \widehat{\widetilde{\bm\lambda}}  
  \Big[\varphi(\mathbf{X}, t_0) - \widehat{\varphi}^{t_0} \Big] \, d(\mathcal{B}\times t_0) = 0\,\qquad \forall \widehat{\widetilde{\bm\lambda}} \in L^2(\mathcal{B}\times t_0) \;.
  \label{eq:weakFlowmapConstraint}
\end{equation}

\begin{figure}[H]
    \centering
    \begin{minipage}{0.4\textwidth}
        \centering
        \includegraphics[width=\linewidth]{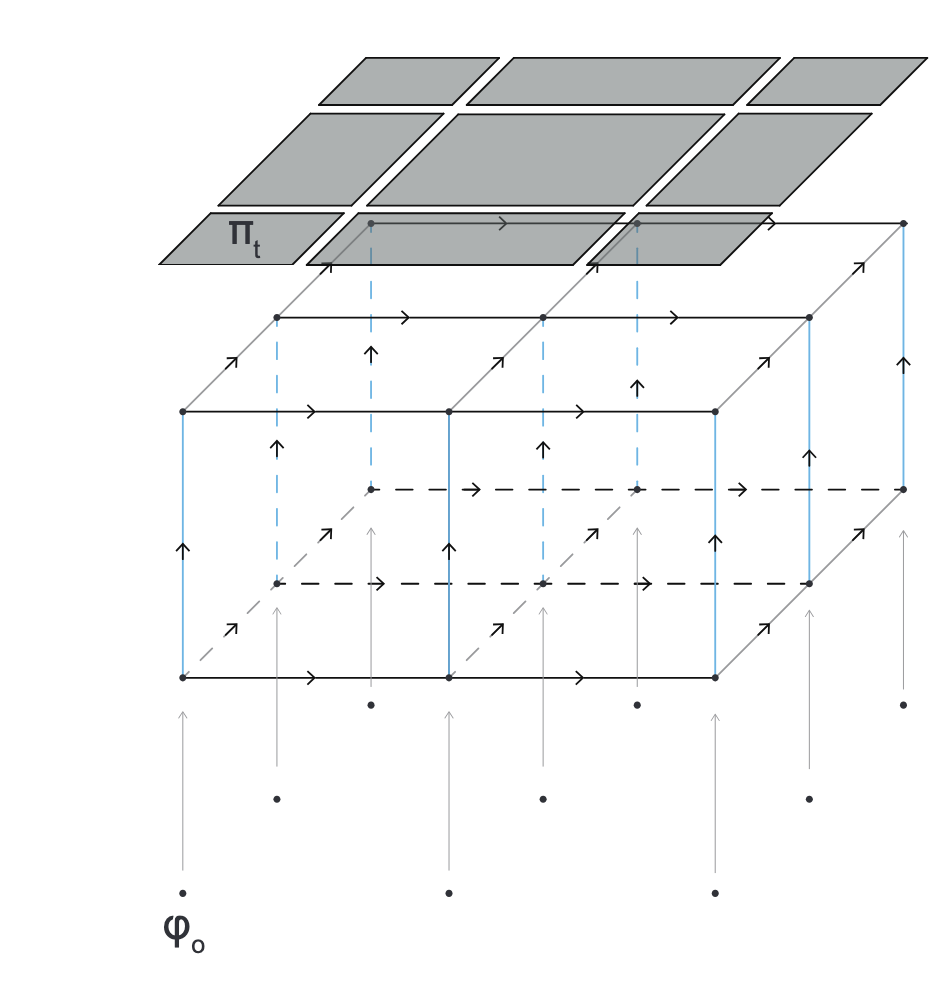}
        \caption{Prescription of initial flow map}
        \label{fig:lm_mom}
    \end{minipage}
    \hspace{2cm}
    \begin{minipage}{0.4\textwidth}
        \centering
        \includegraphics[width=\linewidth]{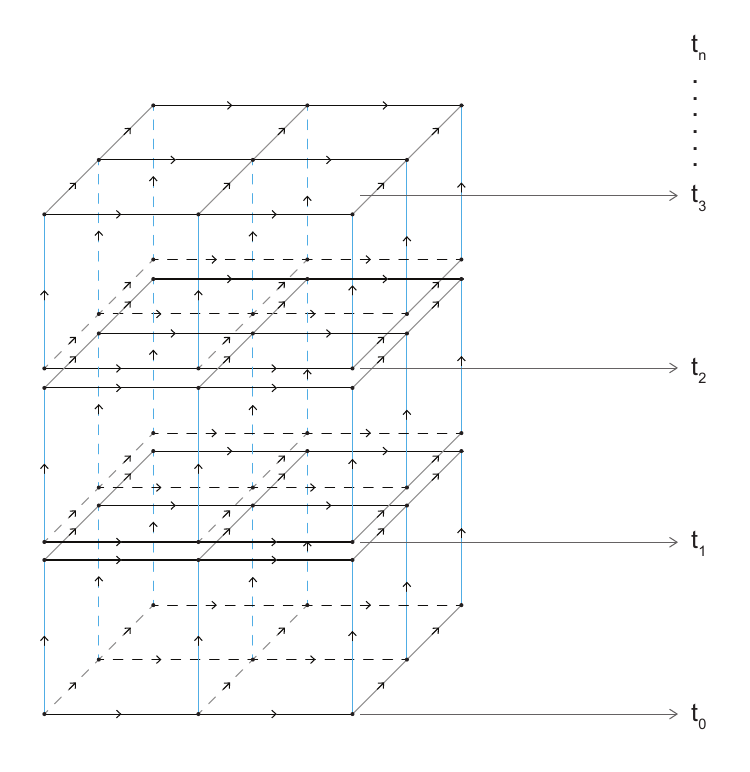}
        \caption{Time stacking}
        \label{fig:time_stacking}
    \end{minipage}
\end{figure}

The simulation proceeds via time stacking, where space-time grids are treated as sequential spectral elements. The computed flow map and momentum (\(\widehat{\bm\pi}^{t^n}\)) from the final time of one element serve as initial conditions for the subsequent element. Finally the equations in \eqref{eq:final_mom_eqn} and \eqref{eq:final_mom_cons_eqn}, upon incorporation of the initial conditions \eqref{eq:weakFlowmapConstraint}, can be assembled into a final matrix form as shown below. \\

\begin{equation}
    \begin{bmatrix}
        \mathbb{M}_{\rho_0}^{-1} & 0 & -\mathbb{E}_t^{1,0} & 0 & 0 & 0   \\[8pt]
        0 & \mathbb{M}_{\rho_o}^{-1} & 0 & -\mathbb{E}_t^{1,0} & 0 & 0  \\[8pt]
        (\mathbb{E}_t^{1,0})^\top & 0 & 0 & -\mathbb{D}_{PK} & -\mathbb{N}^{t^{n+1}}_{\pi} & 0 \\[8pt]
        0 & (\mathbb{E}_t^{1,0})^\top & \; \mathbb{D}_{PK} & 0 & 0 & -\mathbb{N}^{t^{n+1}}_{\pi} \\[8pt] 
        0 & 0 & (\mathbb{N}_{\pi}^{t^{n+1}} )^\top & 0 & 0 & 0  \\[8pt]
        0 & 0 & 0 & (\mathbb{N}_{\pi}^{t^{n+1} })^\top & 0 & 0  \\
    \end{bmatrix}
    \begin{bmatrix}
        \overline{ \widetilde{\bm \pi}}_x  \\[10pt]
        \overline{\widetilde{\bm \pi}}_y \\[10pt]
        \overline{\bm \varphi}_x \\[10pt]
        \overline{\bm \varphi}_y \\[10pt]
        \widehat{ \widetilde{\bm \pi}}_x^{t^{n+1}} \\[10pt]
        \widehat{ \widetilde{\bm \pi}}_y^{t^{n+1}} 
    \end{bmatrix}
    =
    \begin{bmatrix}
        0   \\[10pt]
        0 \\[10pt]
        -\mathbb{N}^{t_o}\widehat{\widetilde{\bm{\pi}}}_x^{t^n}\\[10pt]
        -\mathbb{N}^{t_o}\widehat{\widetilde{\bm{\pi}}}_y^{t^n}\\[10pt]
        \widehat{\varphi}_x^{t_0} \\[10pt]
        \widehat{\varphi}_y^{t_0}
    \end{bmatrix}.
    \label{eq:square_discrete_system}
\end{equation}

Note, here $ \mathbb{D}_{PK}$ is the simplified matrix notation of \eqref{eq:discrete_pk}, \(\mathbb{D}_{PK} = \mathbb{E}_x^{1,0^T} \; \mathbb{M}_{P_w}\;\mathbb{E}^{1,0}_y \  - \;  \; \mathbb{E}_y^{1,0^T} \;  \mathbb{M}_{P_w}\;   \mathbb{E}_x^{1,0}  \). This final, square matrix system \eqref{eq:square_discrete_system}, represents the fully discretized problem for a single space-time element. It is this algebraic system of equations that is solved computationally to find the unknown degrees of freedom.

\section{Conservation Properties of the Formulation}
\label{sec:conservation_properties}

We will use the weak form of the momentum balance, \eqref{eq:final_weak_formulation}, to demonstrate that the formulation inherently conserves key physical quantities. This is shown by selecting specific test functions \(\bm{\phi}\) that correspond to rigid-body motions  in the weak for of the momentum equation.

\subsection{Conservation of Linear Momentum}

\begin{proof}[Proof]
We select a test function representing a rigid-body translation, \(\bm{\phi}(\mathbf{x},t) = \mathbf{c}\), where \(\mathbf{c}\) is an arbitrary constant vector. Substituting this into the weak form \eqref{eq:final_weak_formulation} yields:

\[
     \int\limits_{\mathcal{B}\times T } \widetilde{\bm\pi}\frac{\partial \mathbf{c}}{\partial t} \, d(\mathcal{B}\times T) \;
    - \int\limits_{\mathcal{B}\times \partial T } (\widetilde{\bm\pi} \cdot \bm{n}_t)\;\mathbf{c} \, d(\mathcal{B}\times \partial T) \;
    + \int\limits_{\mathcal{B}\times T } \left ( (P_{W}-P_{ext})J\bm{F}^{-T} \right )\; : \nabla \mathbf{c}  \, d(\mathcal{B}\times T) \;  = \;  0.
\]

Since \(\mathbf{c}\) is constant, its gradient (\(\nabla \mathbf{c}=\mathbf{0}\)) and time derivative (\(\partial \mathbf{c}/\partial t = \mathbf{0}\)) both vanish. The first and third terms vanish, leaving only the boundary term, which simplifies as:

\[
    \int\limits_{\mathcal{B}\times \partial T } (\widetilde{\bm\pi} \cdot \bm{n}_t)\cdot\mathbf{c} \, d(\mathcal{B}\times \partial T) 
      =  \left[ \int\limits_{\mathcal{B}} \mathbf{c} \cdot \widehat{\widetilde{\bm{\pi}}} \,d\mathcal{B} \right]_{t=0}^{t=T} = 0.
\]

Because \(\mathbf{c}\) is arbitrary, we can choose $\bm{c}$ to be the unit vector in the $x$- or $y$-direction, which implies that the total linear momentum \(\int_{\mathcal{B}} \widetilde{\pi}_x \,d\mathcal{B}\) and \(\int_{\mathcal{B}} \widetilde{\pi}_y \,d\mathcal{B}\) is conserved:

\[
    \left( \int\limits_{\mathcal{B}} {\widetilde{{\pi}}_x} \;d\mathcal{B} \right)_{t=0} = \left( \int\limits_{\mathcal{B}} {\widetilde{{\pi}}_x} \;d\mathcal{B} \right)_{t=T} \qquad \mbox{and} \qquad \left( \int\limits_{\mathcal{B}} {\widetilde{{\pi}}_y} \;d\mathcal{B} \right)_{t=0} = \left( \int\limits_{\mathcal{B}} {\widetilde{{\pi}}_y} \;d\mathcal{B} \right)_{t=T} \;.
\]
Here we once again assumed that the external pressure is constant. If this is not the case, we need to use \eqref{eq:final_weak_formulation} and we see that the change in total momentum is equal to the work done by the extrenal pressure over the boundary $\partial \mathcal{B} \times T$.
\end{proof}

\subsection{Conservation of Angular Momentum}

\begin{proof}[Proof]
We select a test function corresponding to a rotation of the flow map \({\bm\phi} = [\varphi_y, -\varphi_x]^T\). Substituting this into the weak form \eqref{eq:weak_form_start} requires us to check each term.

\begin{itemize}
    \item Stress Term: The inner product of the stress \(\widetilde{\bm{P}} = (P_W - P_{ext}) J \mathbf{F}^{-T}\) and the test gradient \(\nabla {\phi}\) vanishes due to the skew-symmetric structure:
    \[
    \int\limits_{\mathcal{B}\times T} (P_W-P_{ext}) \;\left(\;\cancel{F_{yy}\;F_{yx}}\; -\; \cancel{F_{yx}\;F_{yy}} \;+\; \cancel{F_{xy}\;F_{xx}}\; -\; \cancel{F_{xx}\;F_{xy}}\; \right) \ d(\mathcal{B}\times T) = 0.
    \]
    \item Momentum Term: Using the constitutive relation \(\widetilde{\bm{\pi}} = \mathbb{M}_{\rho_0}^{-1} \dot{\bm{\varphi}}\), the momentum term \(\int (\pi_x \dot{\varphi}_y - \pi_y \dot{\varphi}_x) \, dt\) also vanishes due to the symmetry of the weighted mass matrix:
    \[
     {\dot{\overline{\varphi}}_y}^T \mathbb{M}_{\rho_0}^{-1} \dot{\overline{\varphi}}_x - {\dot{\overline{\varphi}}_x}^T \mathbb{M}_{\rho_0}^{-1} \dot{\overline{\varphi}}_y = 0.
    \]
\end{itemize}

With the stress and internal momentum terms eliminated, only the boundary term remains:
\[
 \int\limits_{\mathcal{B} \times \partial T} (\widetilde{\bm{\pi}} \cdot \bm{n}_t) \cdot {\varphi} \, d(\partial \Omega) = \left[ \int\limits_{B} (\widetilde{\pi}_x \varphi_y - \widetilde{\pi}_y \varphi_x) \, d\mathcal{B} \right]_{t=0}^{t=T} = 0.
\]
Therefore, total angular momentum is conserved:
\begin{equation}
    \left( \int\limits_{B} (\widetilde{\pi}_x \varphi_y - \widetilde{\pi}_y \varphi_x) \, d\mathcal{B} \right)_{t=0} = \left( \int\limits_{B} (\widetilde{\pi}_x \varphi_y - \widetilde{\pi}_y \varphi_x) \, d\mathcal{B} \right)_{t=T}.
\end{equation}
\end{proof}

\subsection{Conservation of Energy}

\begin{proof}[Proof]

To establish energy conservation, we utilize once again the weak form \eqref{eq:final_weak_formulation} and construct a specific test field $\bm{\phi}$ based on the velocity field $\dot{\bm{\varphi}}$. Although $\dot{\varphi}$ is expanded along temporal edges (representing time derivatives of degree $(N-1)$), we can exactly represent it in the nodal test space (degree $N$, using the fact that every polynomial of degree $(N-1)$ can be represented in a space spanned by polynomials of degree $N$). We do this by sampling the velocity field at the primal nodes and assigning these values as the degrees of freedom for $\phi$.

\begin{equation}
    \int\limits_{\mathcal{B}\times T } \widetilde{\bm\pi}\;\frac{\partial\dot{\bm{\varphi}}}{\partial t} \, d(\mathcal{B}\times T) \ - \int\limits_{\mathcal{B}\times \partial T } (\widetilde{\bm\pi} \cdot \bm{n}_t)\;\dot{\bm{\varphi}} \, d(\mathcal{B}\times \partial T) \
    + \int\limits_{\mathcal{B}\times T } \left ( (P_{W}-P_{ext) }J\bm{F}^{-T} \right )\; :\nabla\dot{\bm{\varphi}}\; \, d(\mathcal{B}\times T) = 0\;\qquad \forall \dot \varphi \in H^1(\mathcal{B}\times T) \;.
\end{equation}



\paragraph{Kinetic Energy:}  The first two terms can be combined to give

\[ \int\limits_{\mathcal{B}\times T } \widetilde{\bm\pi}\;\frac{\partial\dot{\bm{\varphi}}}{\partial t} \, d(\mathcal{B}\times T) - \int\limits_{\mathcal{B}\times \partial T } (\widetilde{\bm\pi} \cdot \bm{n}_t)\;\dot{\bm{\varphi}} \, d(\mathcal{B}\times \partial T) = \int\limits_{\mathcal{B}\times T } \frac{\partial \widetilde{\bm\pi}}{\partial t} \;\dot{\bm{\varphi}} \, d(\mathcal{B}\times T) \;.\]

Now using \(\bm{\pi} = \rho_0 \dot{\varphi}\) this can be written as
 
 \[
  \int\limits_{\mathcal{B}\times T } \frac{\partial \widetilde{\bm\pi}}{\partial t} \;\dot{\bm{\varphi}} \, d(\mathcal{B}\times T)  =
  \int\limits_{\mathcal{B}\times T}
     \rho_0\,\frac{\partial}{\partial t}\Bigl(\frac{\partial \varphi}{\partial t}
    \,\cdot\,\frac{\partial \varphi}{\partial t}\Bigr)
  \;d(\mathcal{B}\times T) = \int\limits_{\mathcal{B} \times T} \frac{\partial}{\partial t} \left ( \frac{1}{2} \rho_0 \frac{\partial \bm{\varphi}}{\partial t} \cdot \frac{\partial \bm{\varphi}}{\partial t}  \right ) \ \mathrm{d}(\mathcal{B}\times T) \;.
  \
\]




Hence the first two terms in the weak formulation denote  the rate of change of kinetic energy within the space-time element,

 \begin{equation}
 \int\limits_{\mathcal{B}\times T } \widetilde{\bm\pi}\;\frac{\partial\dot{\bm{\varphi}}}{\partial t} \, d(\mathcal{B}\times T) - \int\limits_{\mathcal{B}\times \partial T } (\widetilde{\bm\pi} \cdot \bm{n}_t)\;\dot{\bm{\varphi}} \, d(\mathcal{B}\times \partial T)
 = \int\limits_{\mathcal{B}\times T} \frac{\partial}{\partial t} \left( \underbrace{\frac{1}{2} \rho_0 \,\bigl\|\dot{\varphi}\bigr\|^2}_{\text{Kinetic Energy}} \right) \;d(\mathcal{B}\times T).
\label{eq:Kinetic_energy_change}
\end{equation}

\paragraph{Internal Energy:} The second term, becomes the power density. For an inviscid barotropic assumption, the first Piola-Kirchoff stress tensor can be written as the the derivative of the internal energy with respect to deformation gradient, this gives \cite{Reddy}.
    
    \[ \left ({{P}_W}-P_{ext} \right ) J \bm{F}^{-\top} : \dot{\bm{F}} \;=\; \frac{\partial W}{\partial \bm F} :\dot{\bm{F}} \;=\;\frac{\partial }{\partial t} W(\bm{F}) \; . \]

Hence, the integral, can be rewritten as the rate of change of internal energy,

    \begin{equation}
     \int\limits_{\mathcal{B}\times T }  \left ({{P}_W}-P_{ext} \right ) J \bm{F}^{-\top} : \dot{\bm{F}}\; \, d(\mathcal{B}\times T)  = \int\limits_{\mathcal{B}\times T }  \;\frac{\partial }{\partial t} W(\bm{F})  \;d(\mathcal{B}\times T)
     \label{eq:internal_energy_change}
    \end{equation}

Combining both the equations \eqref{eq:Kinetic_energy_change} and \eqref{eq:internal_energy_change} gives the total energy conservation statement,

\begin{equation}
    \int\limits_{\mathcal{B}\times T}  
    \frac{\partial}{\partial t} \bigg(  \;
        \frac{1}{2} \rho_0 \bigl\lVert \dot{\varphi} \bigr\rVert^2 
        +  
        W(\bm{F})
    \bigg)  
    \, d(\mathcal{B}\times T) = 0.
\end{equation}
This demonstrates that the sum of the kinetic and internal (potential) energies, thus the total energy remains constant in time.

\end{proof}

\section{Numerical Experiments and Results}\label{sec:NumericalResults}

This chapter presents computational experiments to test the discretization's ability to handle pressure boundary conditions and to verify the conservation of key quantities (mass, momentum, and energy). Two low-Mach test cases are examined: an expansion case ($\alpha = 0.85$) and a compression case ($\alpha = 1.15$). Since benchmark solutions are unavailable, validation proceeds by examining the solution’s convergence as the polynomial order increases from $N=2$ to $N=5$.

\subsection{Experimental Setup}

A $1 \times 1\,\mathrm{m}$ fluid domain is used, with \(\rho = 1.25\,\mathrm{kg/m^3}\), \(\gamma = 1.4\), and \(P_{\text{ref}} = 1\,\mathrm{Pa}\). The boundary pressure is set by \(P_{\text{env}} = \alpha P_{\text{ref}}\). The spatial polynomial order is varied (\(N = 2, 3, 4, 5\)), while the temporal order is \(N_{t} = 1\). The simulation runs for $7.0\,\mathrm{s}$ with \(\Delta t = 0.01\,\mathrm{s}\).

A Picard (fixed-point) iteration is used at each time step to solve the non-linear system, reassembling the pressure mass matrix \(\mathbb{M}_{P_{w}}\) until changes in the deformation (Jacobian $J$ and trace) fall below a tolerance of \(tol = 10^{-12}\). Key parameters are summarized in Table \ref{tab:sim-params}. \\

\begin{table}[h!]
\centering
\renewcommand{\arraystretch}{1.5}
\begin{tabular}{|l|c|l|}
\hline
\textbf{Parameter} & \textbf{Symbol} & \textbf{Value} \\ 
\hline
Fluid density & \(\rho\) & \(1.25\,\mathrm{kg\,m^{-3}}\) \\
Specific heat ratio & \(\gamma\) & \(1.4\) \\
Reference pressure & \(P_{\text{ref}}\) & \(1\,\mathrm{Pa}\) \\
Boundary pressure ratio & \(\alpha\) & \(0.85,\,1.15\) \\
Spatial polynomial orders & \(N\) & \(2,\,3,\,4,\,5\) \\
Temporal polynomial order & \(N_{t}\) & \(1\) \\
Convergence tolerance & $ tol $ & $10^{-12}$ \\
Time step size & \(\Delta t\) & \(0.01\,\mathrm{s}\) \\
Total simulated time & $T_{final}$ & $7\ \mathrm{s}$ \\
\hline
\end{tabular}
\vspace{0.2cm}
\caption{Summary of key simulation parameters.}
\label{tab:sim-params}
\end{table}

\subsection{Results}

This section presents the evolution of the flow for the expansion and compression cases, followed by a convergence analysis. Results for the highest polynomial order ($N=5$) are shown as the representative case.

\subsubsection{Expansion Case: $\alpha = 0.85$}

\begin{figure}[H]
    \centering
    \includegraphics[width=\textwidth]{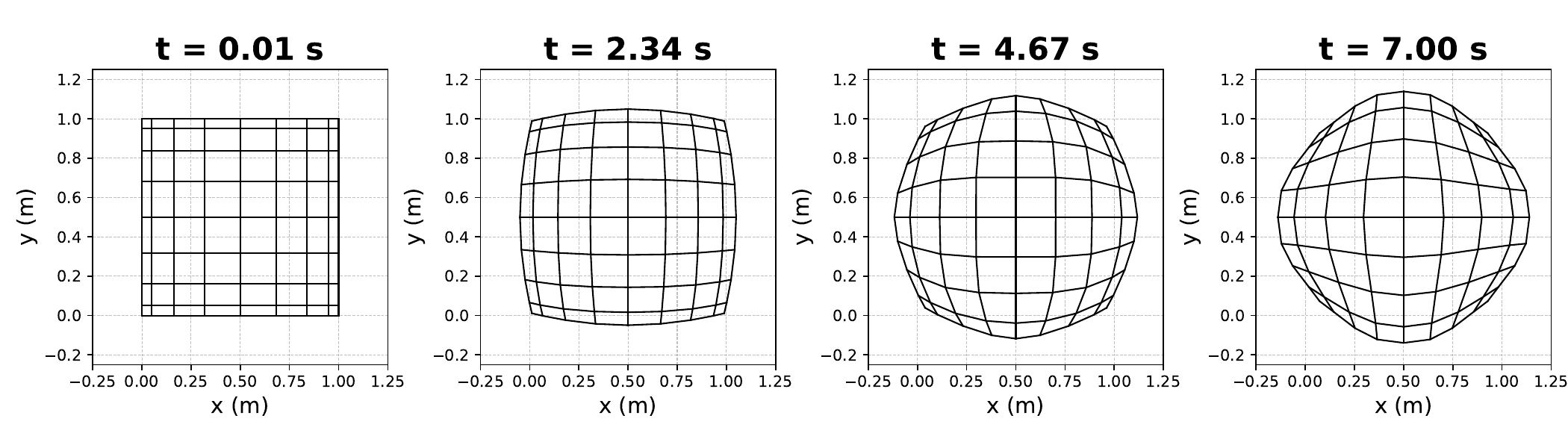}
    \caption{Flow map snapshots for $\alpha=0.85$, $N=5$.}
    \label{fig:flowMap85_N5}
\end{figure}

\begin{figure}[H]
    \centering
    \includegraphics[width=\textwidth]{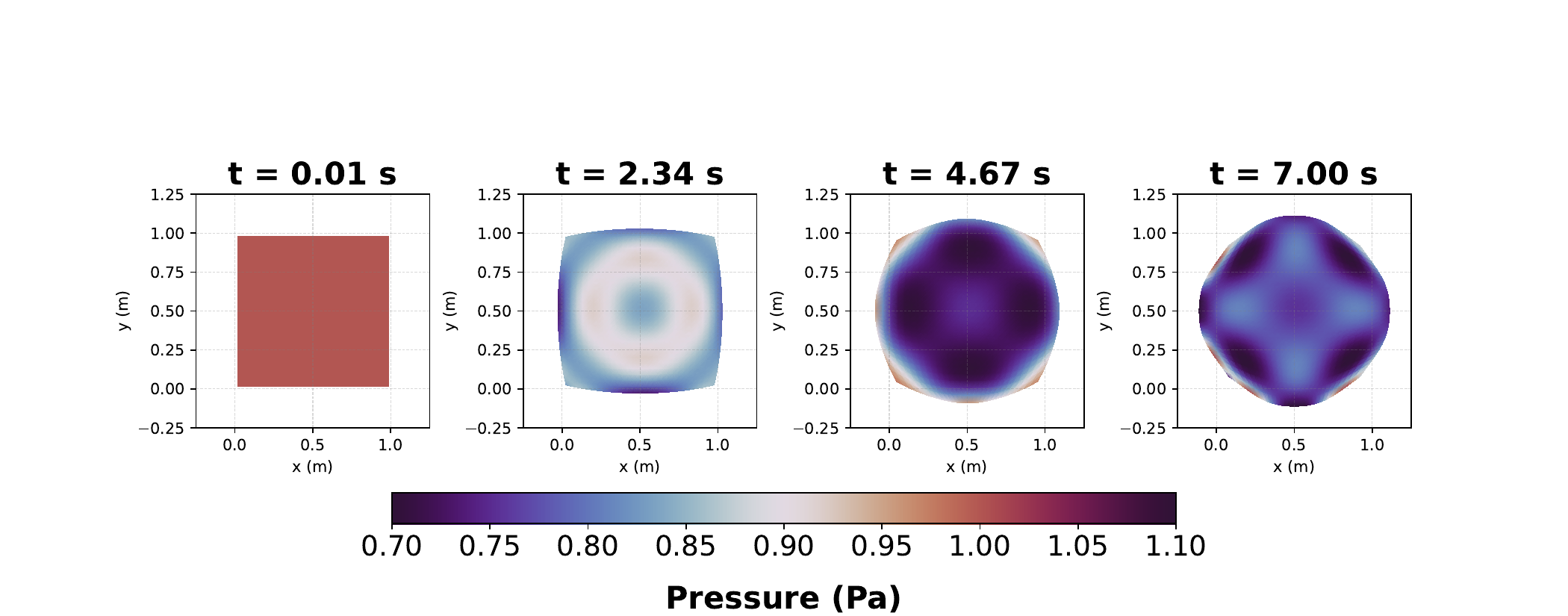}
    \caption{Pressure snapshots (in Pa) for $\alpha=0.85$, $N=5$.}
    \label{fig:pressure85_N5}
\end{figure}

\begin{figure}[H]
    \centering
    \includegraphics[width=\textwidth]{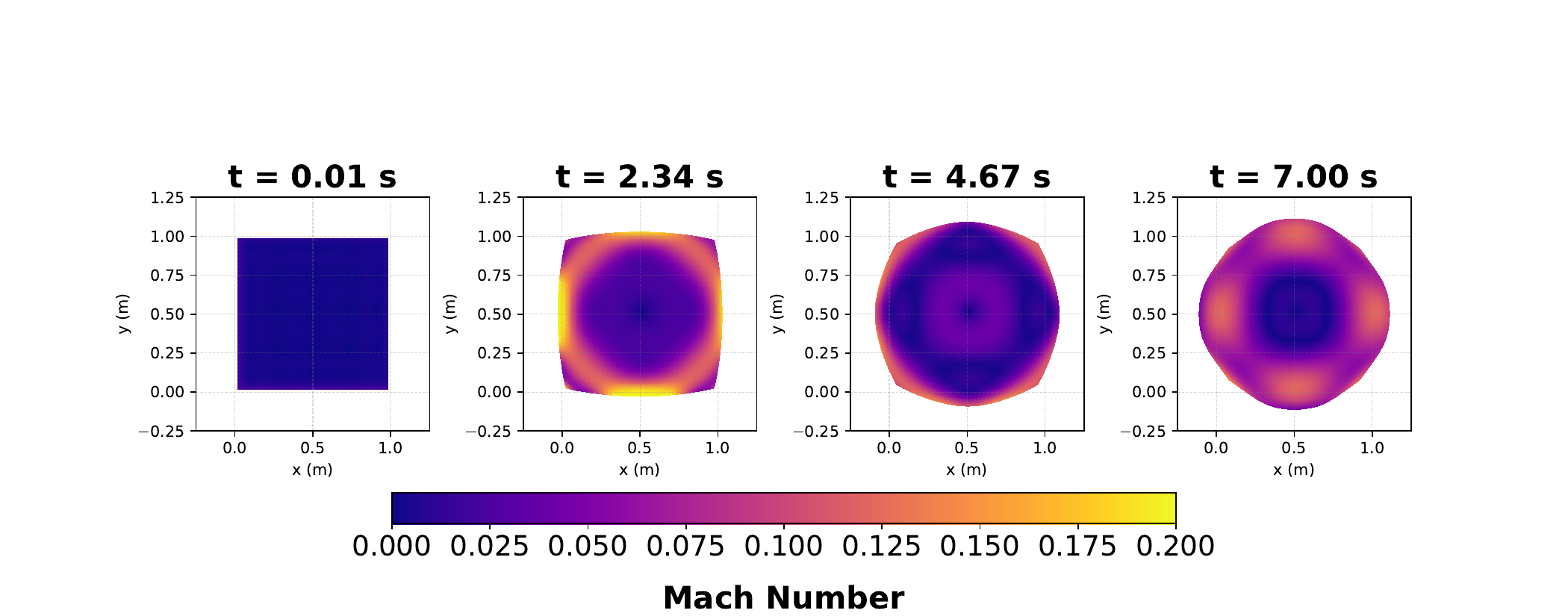}
    \caption{Mach number snapshots for $\alpha=0.85$, $N=5$.}
    \label{fig:mach85_N5}
\end{figure}

\paragraph{Flow Dynamics}
With the boundary pressure set to $0.85\,\mathrm{Pa}$, the fluid expands. Figures ~\ref{fig:flowMap85_N5},~\ref{fig:pressure85_N5},~\ref{fig:mach85_N5} shows the flow map evolution, pressure, and Mach number respectively. Figure~\ref{fig:flowMap85_N5} clearly illustrates the outward deformation. Expansion waves propagate from the boundary and reflect, causing oscillations (Figure~\ref{fig:pressure85_N5}. The Mach number remains low ($\lesssim 0.2$), confirming a nearly incompressible, low-Mach regime (Figure~\ref{fig:mach85_N5}).

\vspace{-1cm}

\vspace{-1cm}

\begin{figure}[H]
    \centering
    \begin{minipage}{\textwidth} 
        \centering
        \begin{subfigure}{0.49\textwidth} 
            \includegraphics[width=\textwidth]{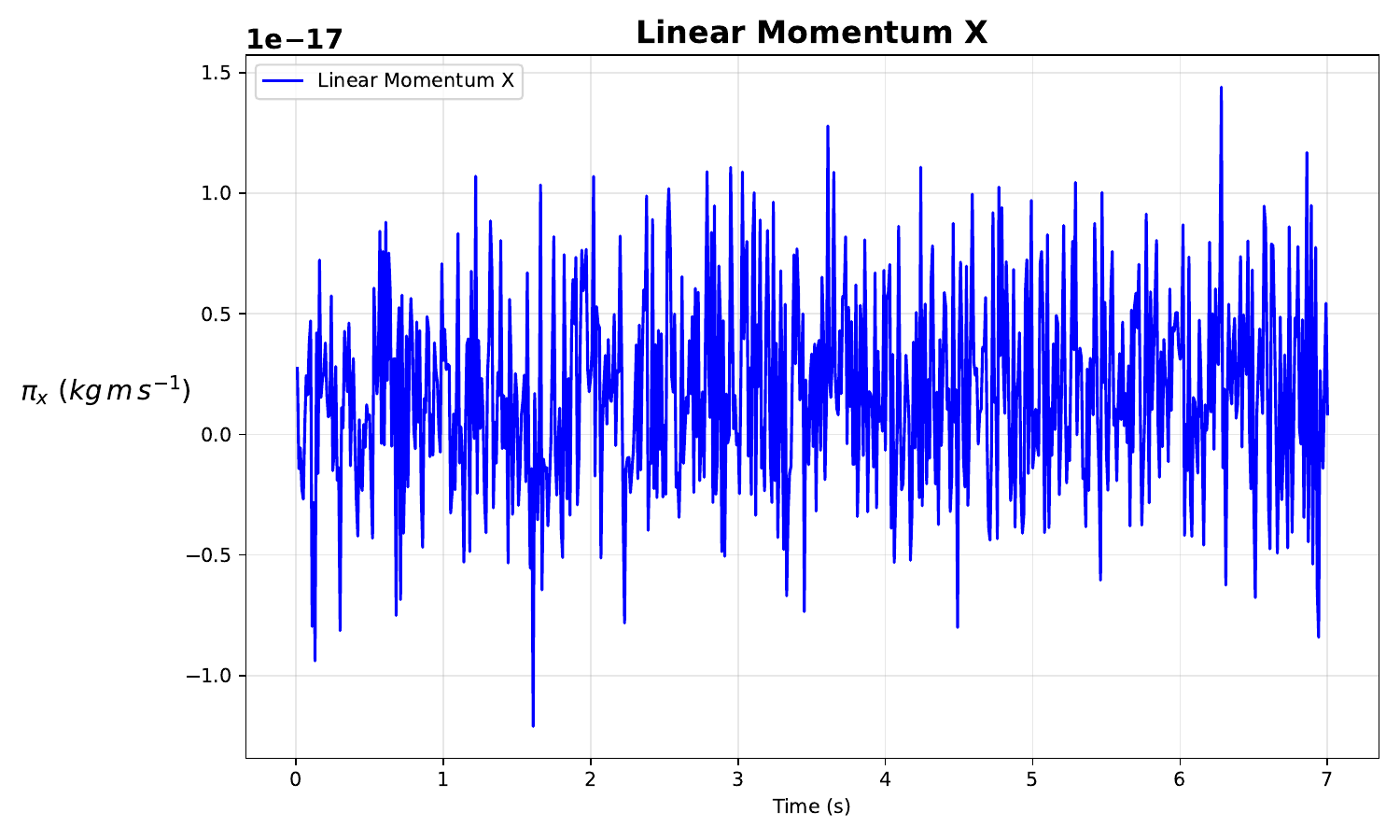}
            \caption{Total linear momentum ($x$).}
        \end{subfigure}
        \hfill %
        \begin{subfigure}{0.49\textwidth}
            \includegraphics[width=\textwidth]{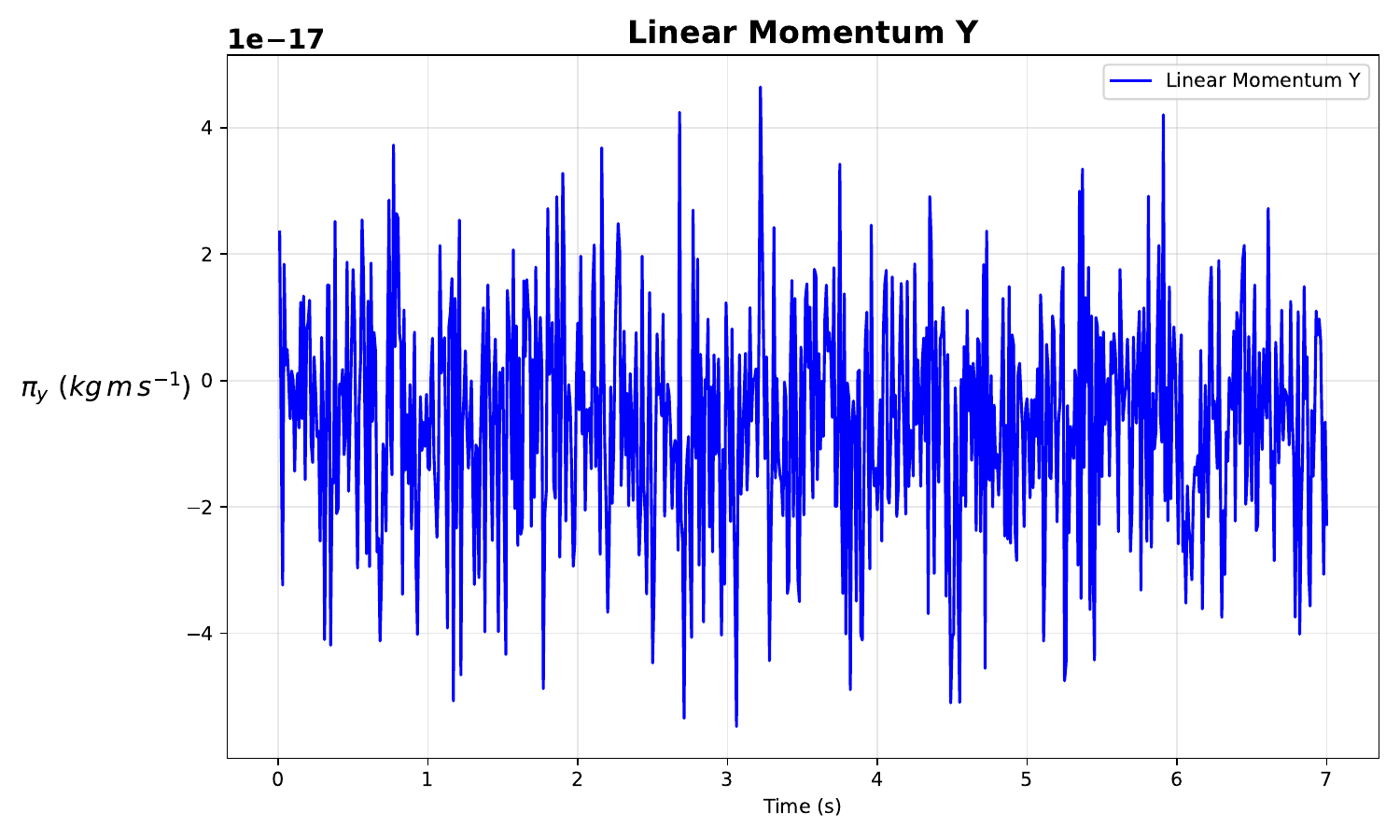}
            \caption{Total linear momentum ($y$).}
        \end{subfigure}
    \end{minipage}
    
    \vspace{7mm} 
  
    \begin{minipage}{\textwidth} 
        \centering
        \begin{subfigure}{0.49\textwidth}
            \includegraphics[width=\textwidth]{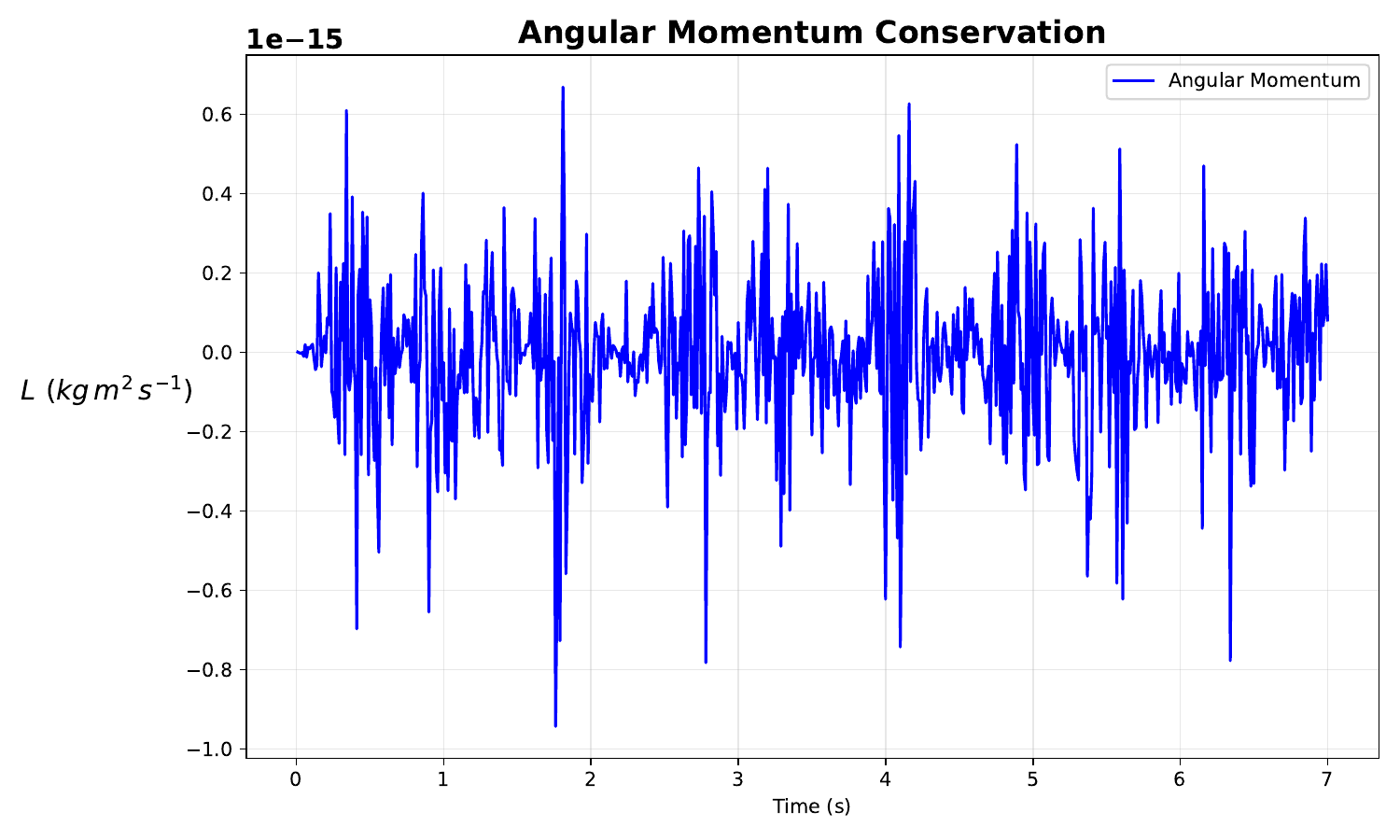}
            \caption{Angular momentum.}
        \end{subfigure}
        \hfill
        \begin{subfigure}{0.49\textwidth}
            \includegraphics[width=\textwidth]{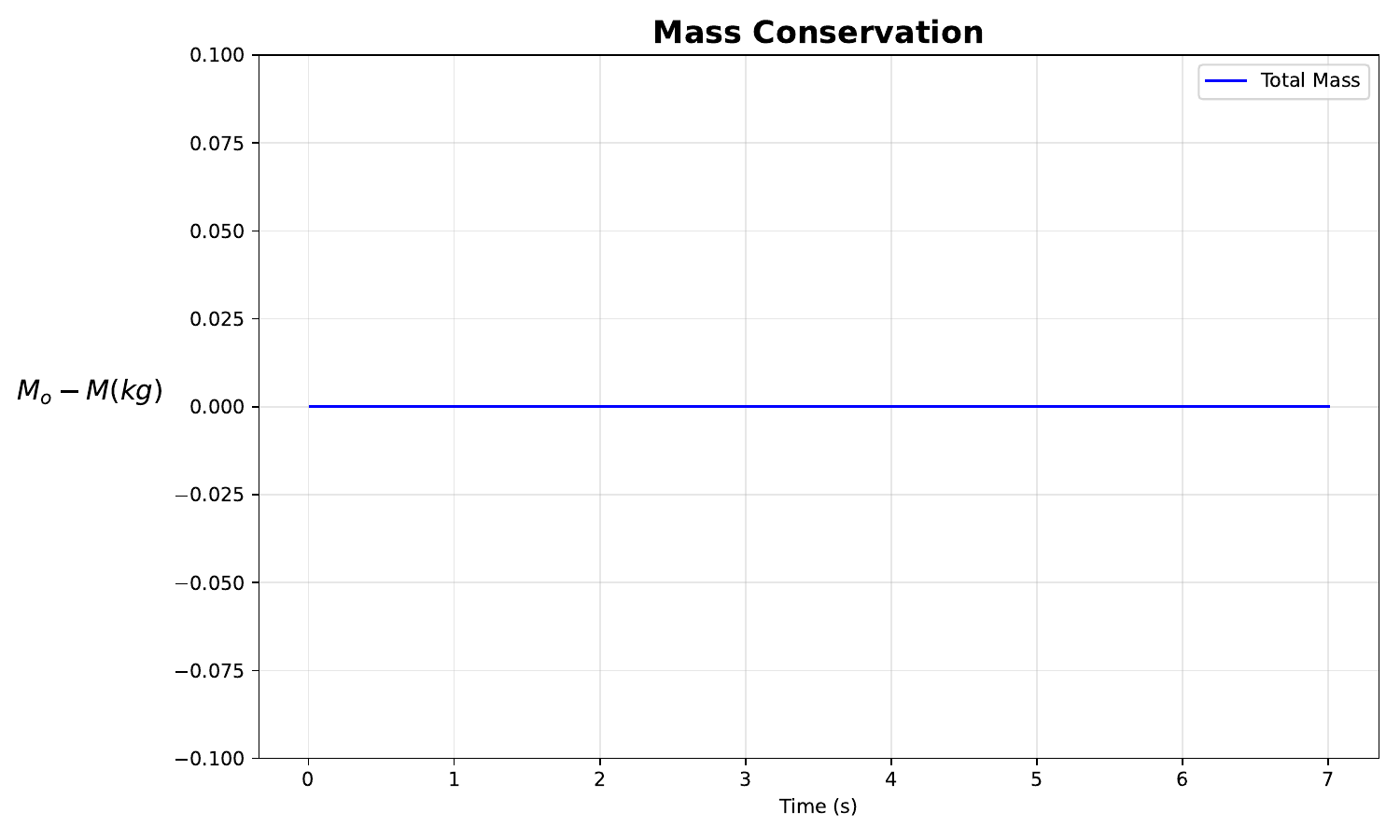}
            \caption{Mass conservation.}
        \end{subfigure}
    \end{minipage}

    \vspace{7mm} 
   
    \begin{subfigure}{0.49\textwidth} 
        \includegraphics[width=\textwidth]{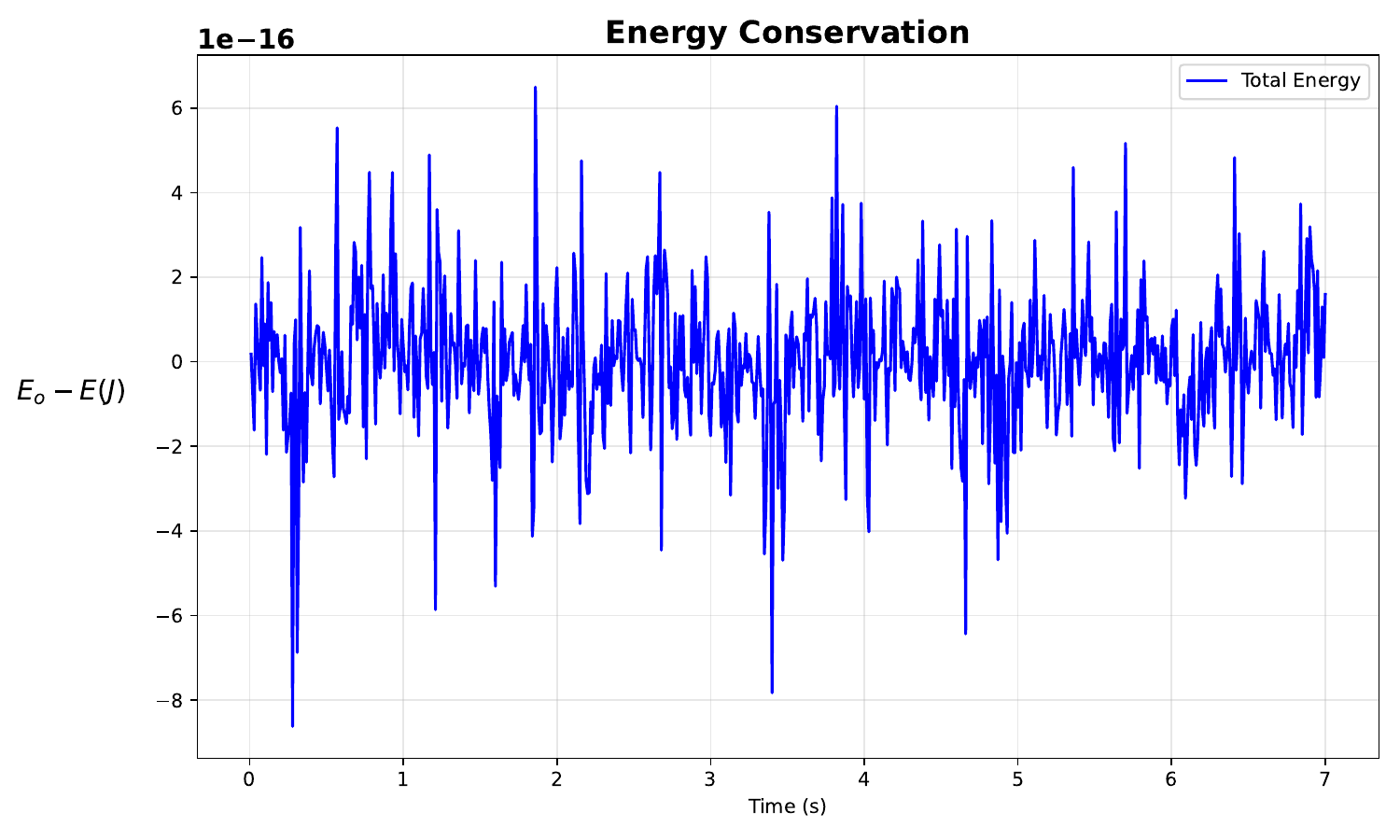}
        \caption{Change in total energy.}
    \end{subfigure}
    
    \caption{Conservation properties for $\alpha=0.85$, $N=5$.}
    \label{fig:conservation_85_N5} 
\end{figure}

\paragraph{Conservation Properties}
The structure-preserving nature of the scheme is demonstrated in Figure~\ref{fig:conservation_85_N5}. Total linear momentum, angular momentum, mass, and energy are all perfectly conserved to machine precision (with fluctuations on the order of $10^{-15}$ to $10^{-17}$) for the entire simulation.

\subsubsection{Compression Case: $\alpha = 1.15$}

\begin{figure}[htbp]
    \centering
    
    \begin{subfigure}{\textwidth}
       
        \makebox[\textwidth][c]{\includegraphics[width=1.1\textwidth]{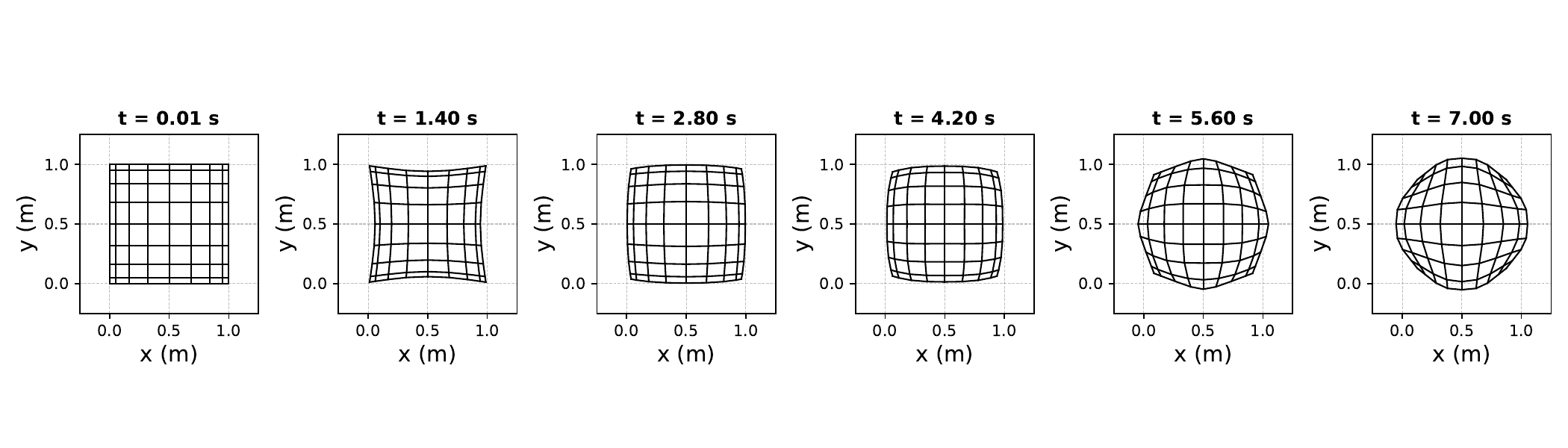}}
        
        \caption{Snapshots of the flow map.}
        \label{fig:flowMap115_N5} 
    \end{subfigure}
    
 \vspace{-1mm}
    
    \begin{subfigure}{\textwidth}
        
        \makebox[\textwidth][c]{\includegraphics[width=1.3\textwidth]{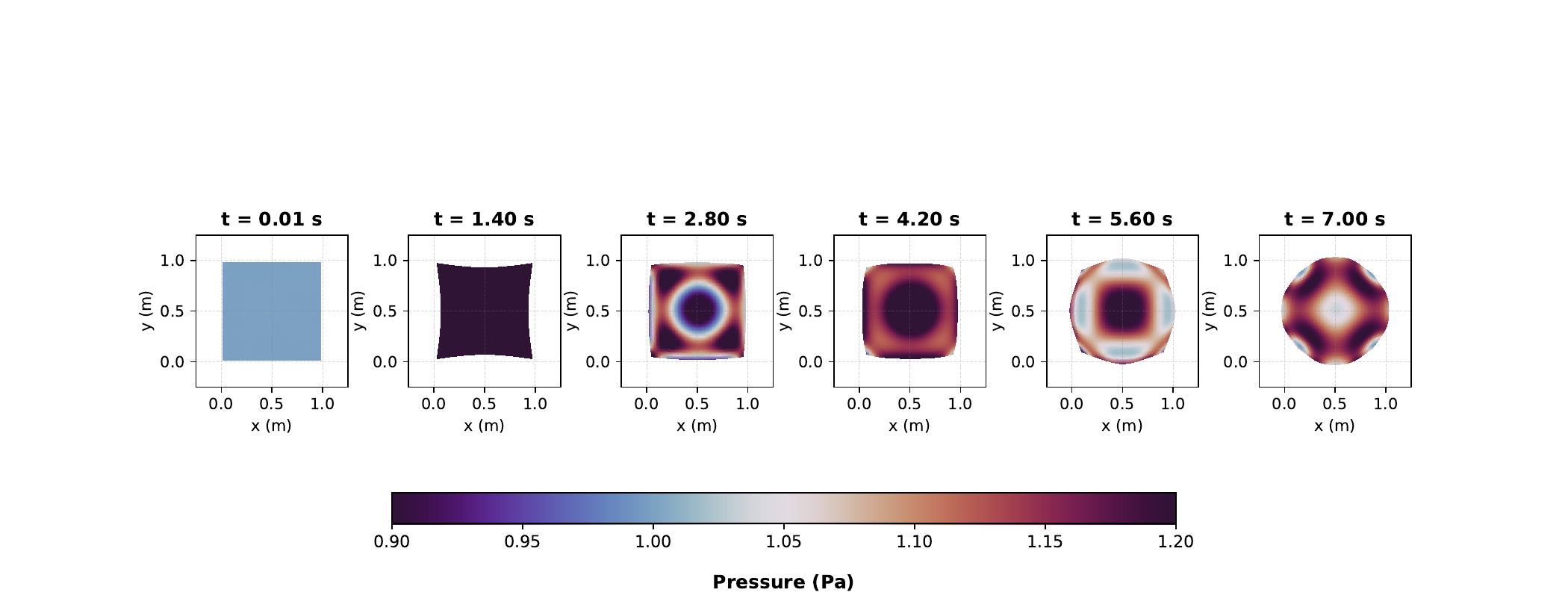}}
        
        \label{fig:pressure115_N5} 
    \end{subfigure}
    
\vspace{-3mm}
    
    \begin{subfigure}{\textwidth}

        \makebox[\textwidth][c]{\includegraphics[width=1.25\textwidth]{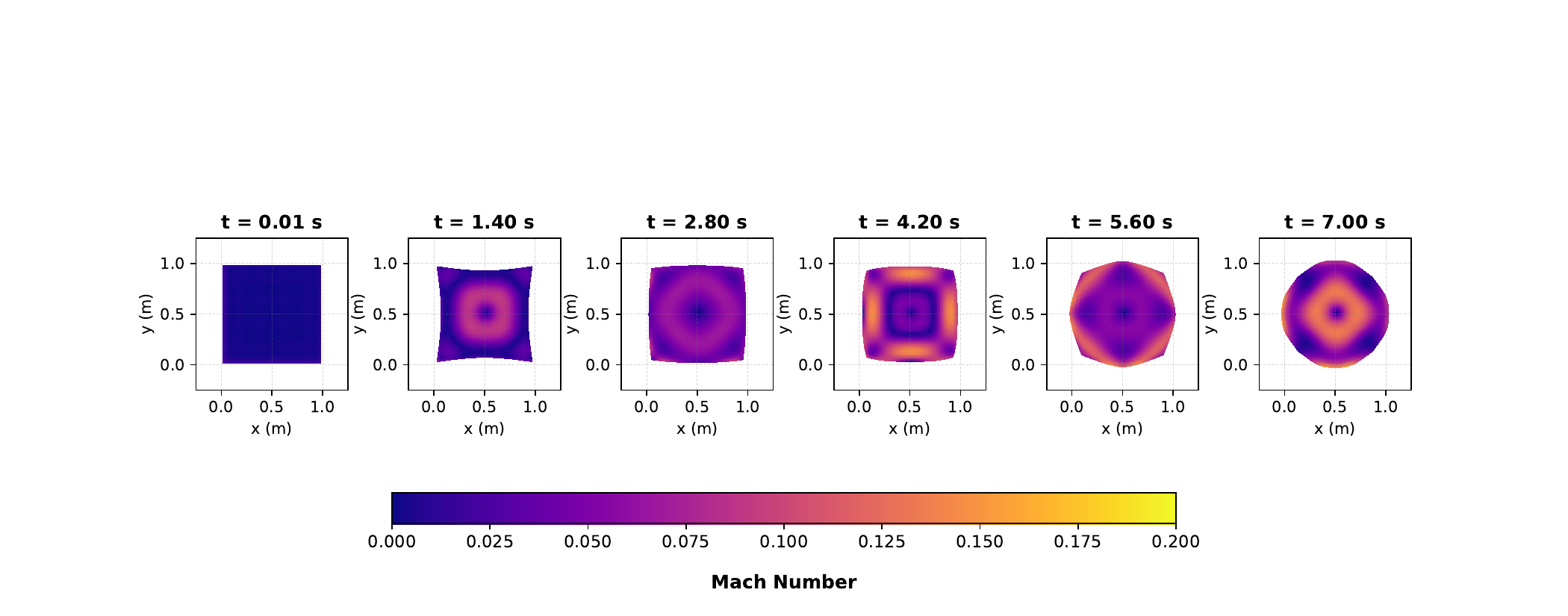}}
        
        \label{fig:mach115_N5}
    \end{subfigure}
    
    \caption{Flow dynamics for $\alpha=1.15$, $N=5$. (a) Flow map, (b) Pressure, (c) Mach number.}
    \label{fig:flow_dynamics_115} 
\end{figure}

\paragraph{Flow Dynamics}
With the boundary pressure set to $1.15\,\mathrm{Pa}$, the fluid is compressed. Figure~\ref{fig:flow_dynamics_115} shows the flow map snapshots, pressure, and Mach number. Figure~\ref{fig:flow_dynamics_115}(a)illustrates the inward deformation. The higher external pressure causes the domain to shrink, with compression waves reflecting internally (Figure~\ref{fig:flow_dynamics_115}(b) and (c)). The flow remains subsonic.

\FloatBarrier

\begin{figure}[htbp]
    \centering
    
    \begin{minipage}{\textwidth}
        \centering
        \begin{subfigure}{0.48\textwidth}
            \includegraphics[width=\textwidth]{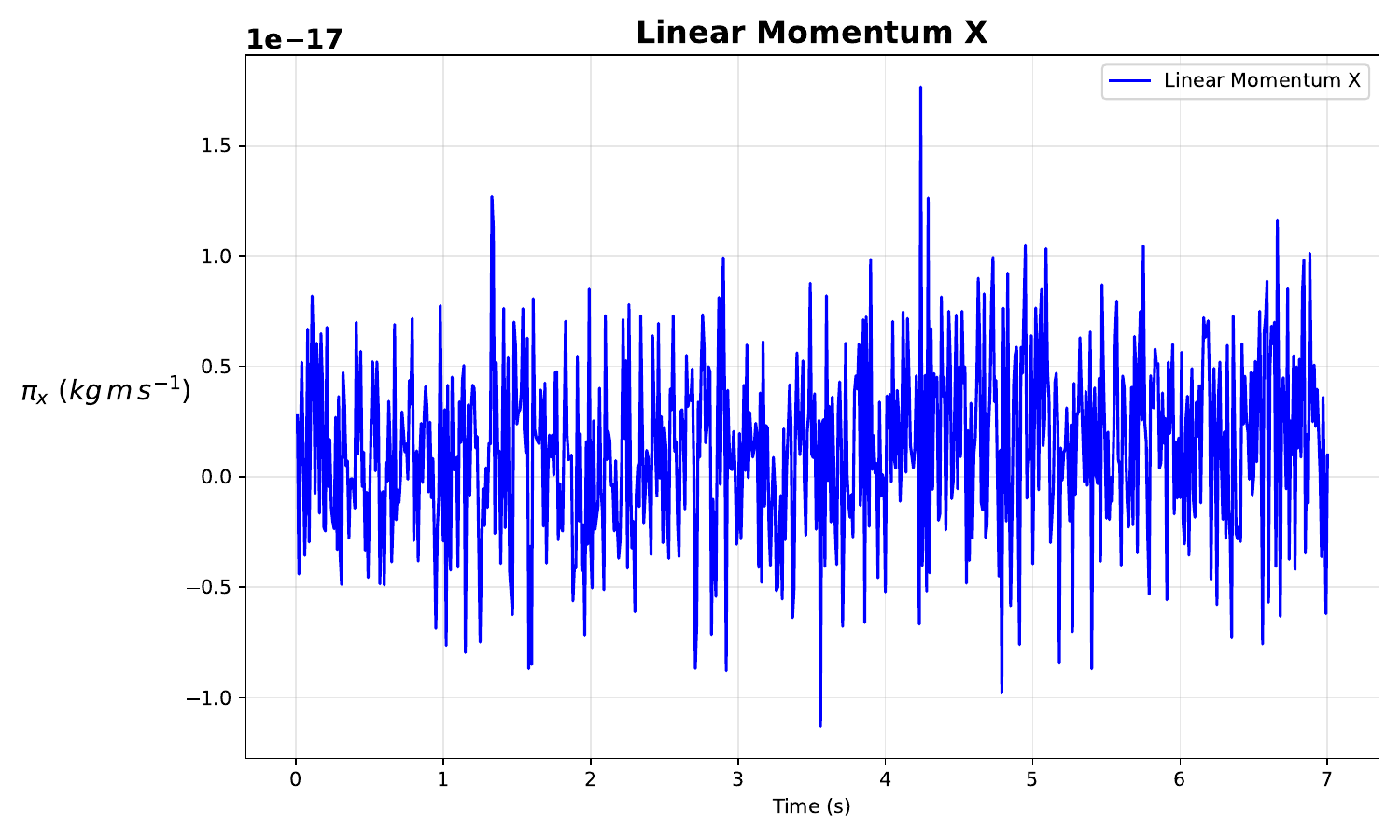}
            \caption{Total linear momentum ($x$).}
        \end{subfigure}
        \hfill 
        \begin{subfigure}{0.48\textwidth}
            \includegraphics[width=\textwidth]{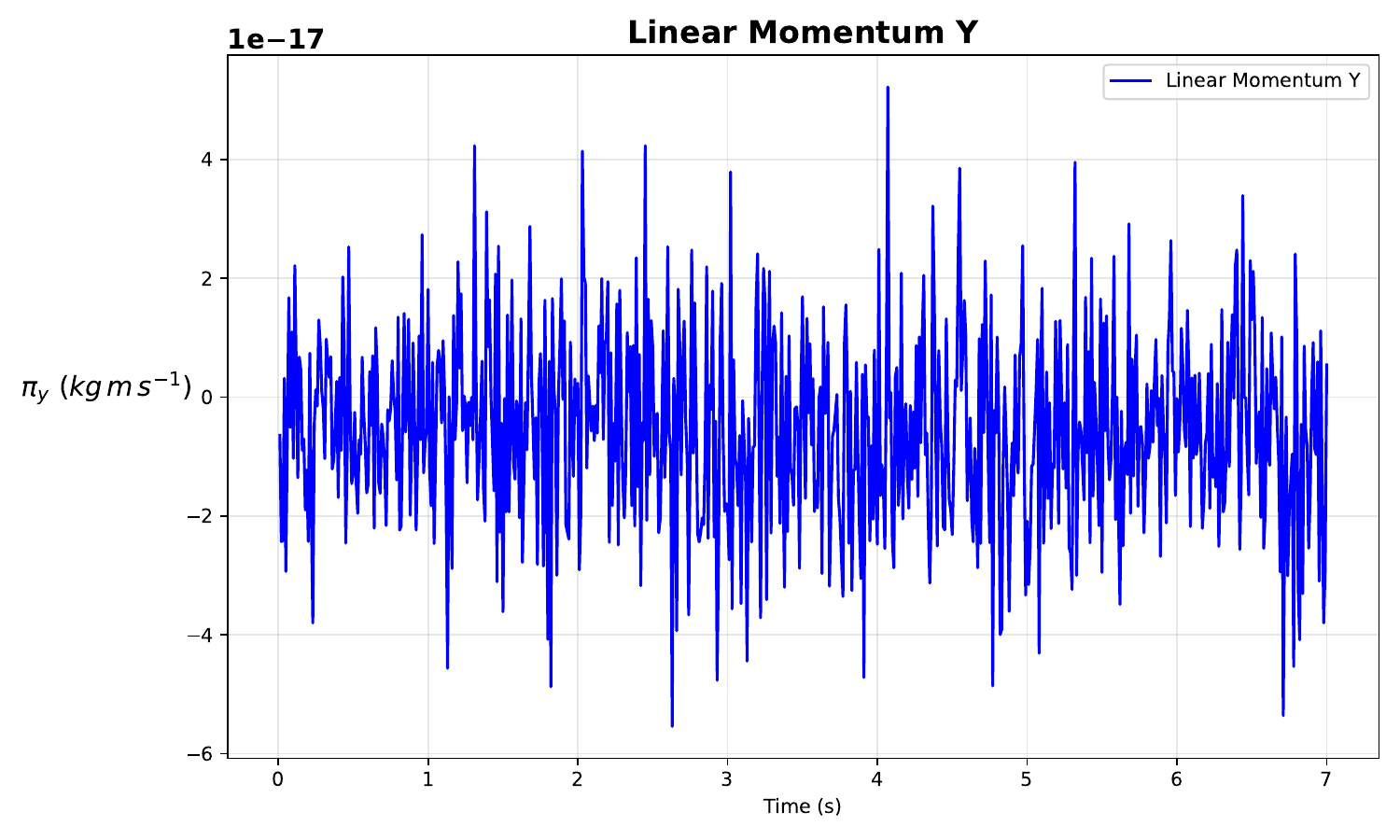}
            \caption{Total linear momentum ($y$).}
        \end{subfigure}
    \end{minipage}
    
    \vspace{7mm} 
    
    \begin{minipage}{\textwidth}
        \centering
        \begin{subfigure}{0.48\textwidth}
            \includegraphics[width=\textwidth]{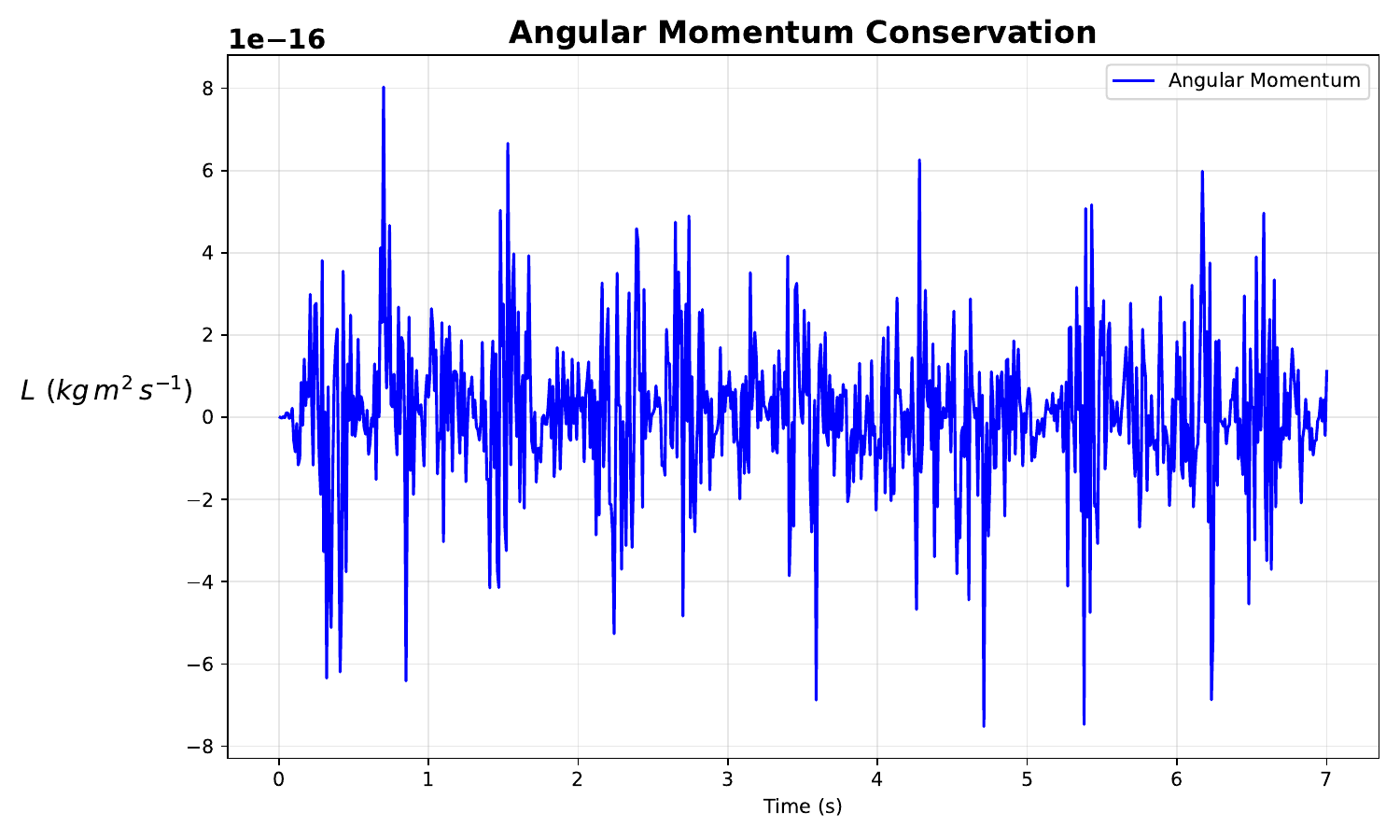}
            \caption{Angular momentum.}
        \end{subfigure}
        \hfill 
        \begin{subfigure}{0.48\textwidth}
            \includegraphics[width=\textwidth]{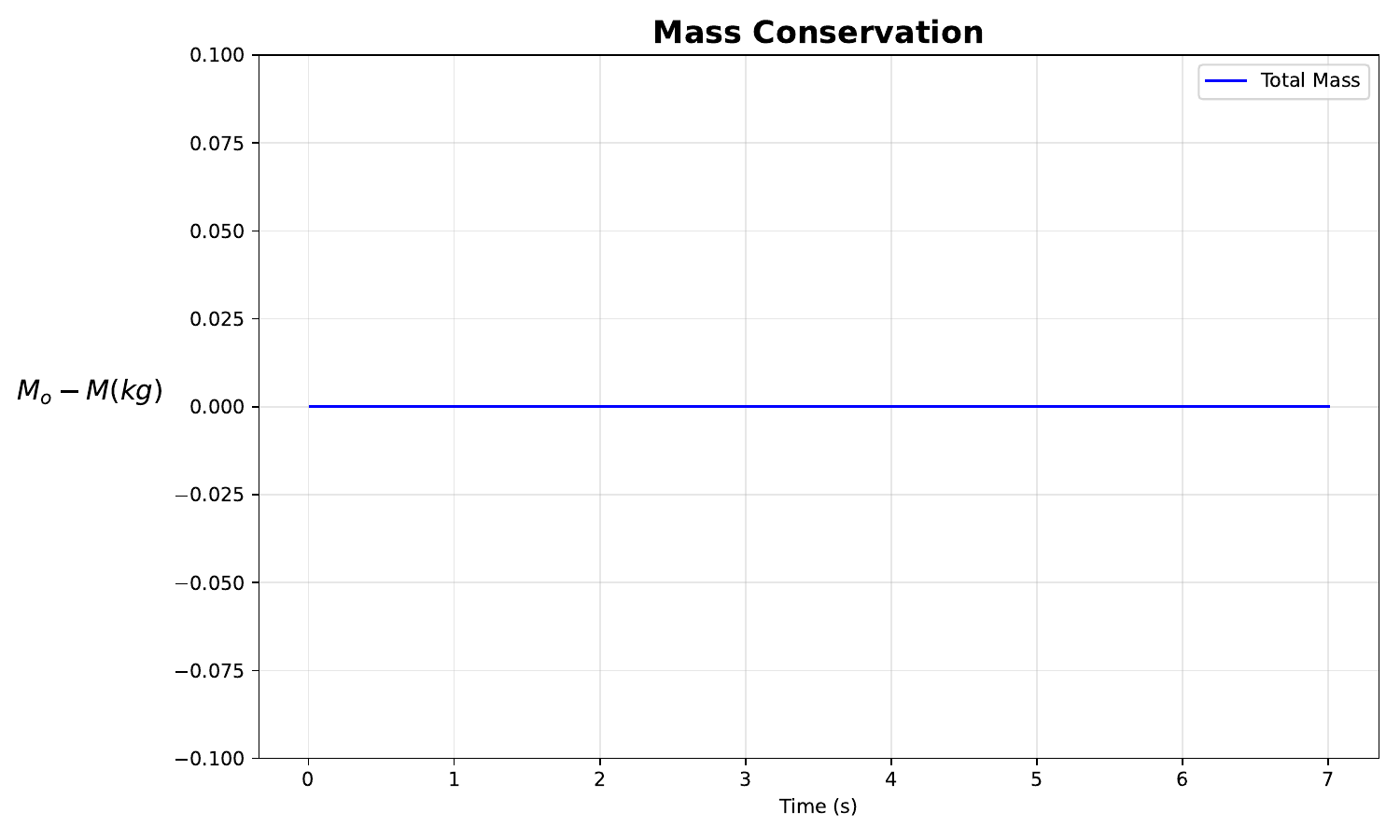}
            \caption{Mass conservation.}
        \end{subfigure}
    \end{minipage}

    \vspace{7mm}

    \begin{subfigure}{0.48\textwidth} 
        \includegraphics[width=\textwidth]{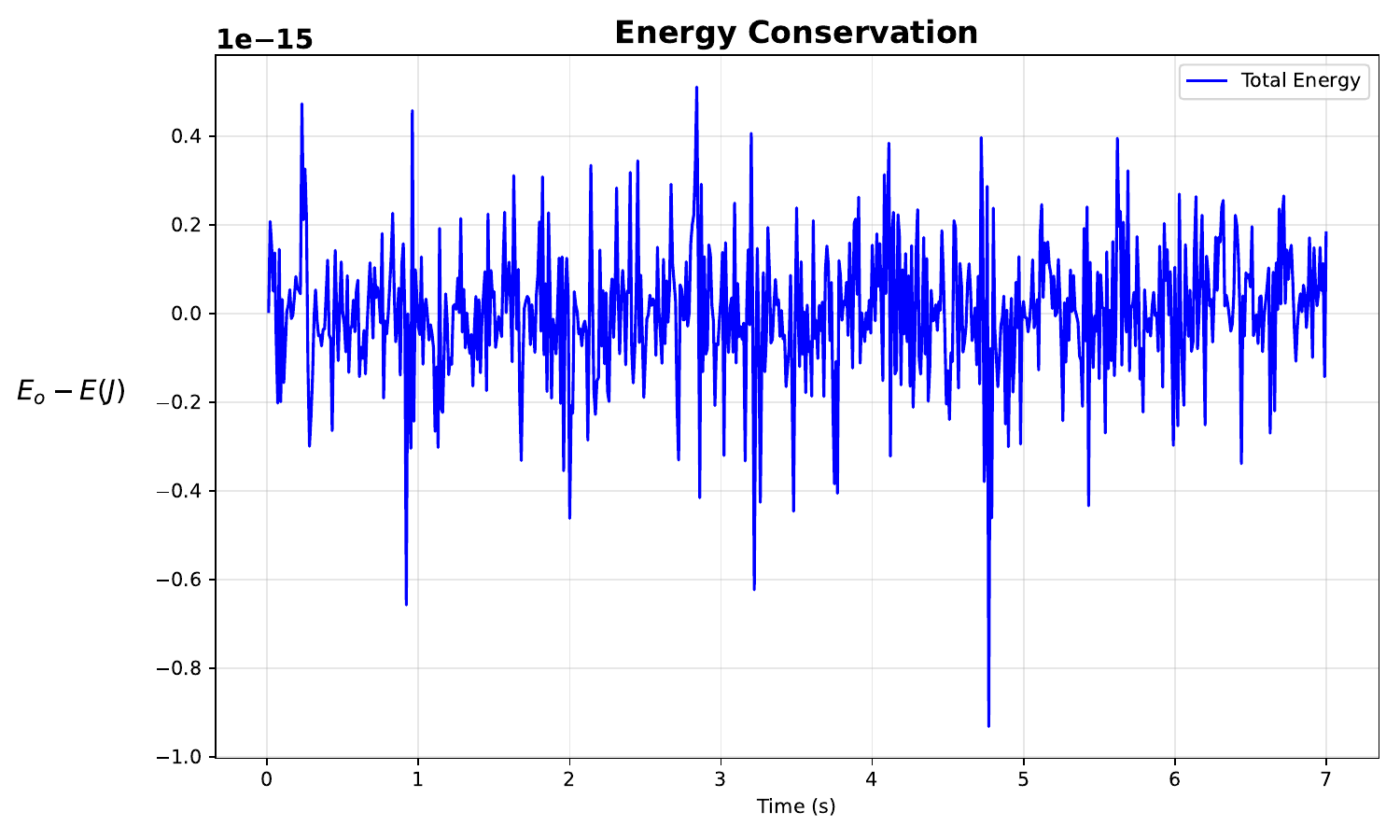}
        \caption{Change in total energy.}
    \end{subfigure}
    
    \caption{Conservation properties for $\alpha=1.15$, $N=5$.}
    \label{fig:conservation_115_N5} 
\end{figure}

\paragraph{Conservation Properties}
As in the expansion case, all physical invariants are conserved to machine precision. Figure~\ref{fig:conservation_115_N5} confirms that mass, linear momentum, angular momentum, and total energy remain constant. \\

\FloatBarrier

\newpage

\subsubsection{Deformation Convergence}

\begin{figure}[htbp] 
    \makebox[\textwidth][l]{\hspace{2cm}\includegraphics[width=0.9\textwidth]{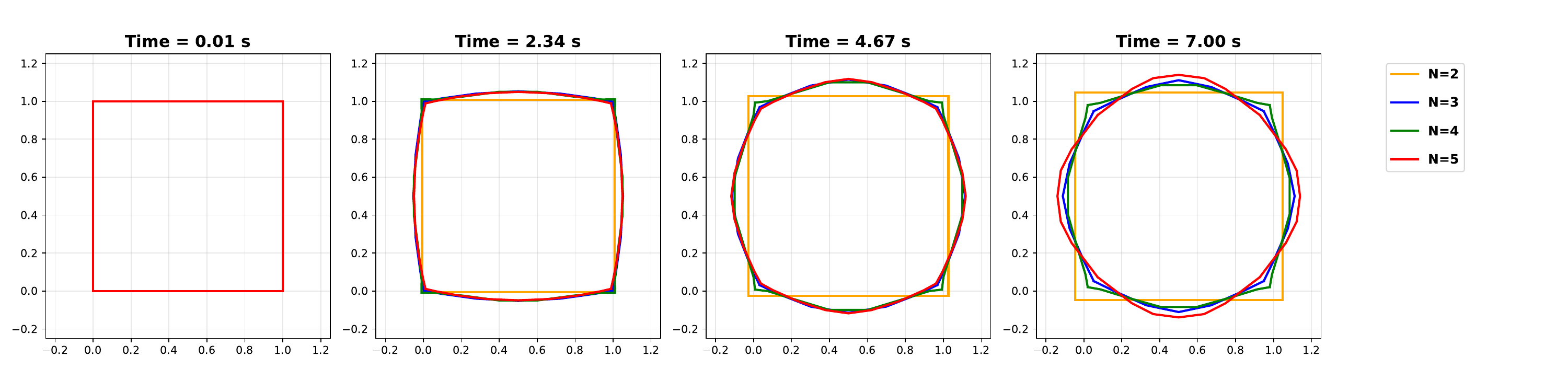}}
    \caption{Convergence of flow map deformation for $\alpha=0.85$ at $t=\{0.01,\,2.34,\,4.67,\,7.00\}\,\mathrm{s}$ for $N=2,3,4,5$.}
    \label{fig:convergence85}
\end{figure}

\begin{figure}[htbp] 
    \makebox[\textwidth][l]{\hspace{2cm}\includegraphics[width=0.9\textwidth]{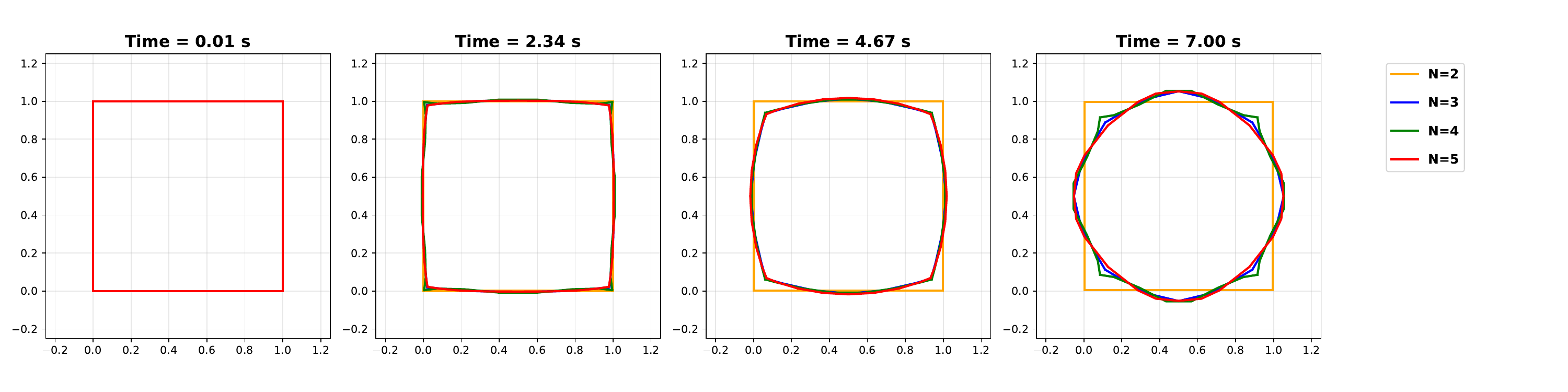}}
    \caption{Convergence of flow map deformation for $\alpha=1.15$ at $t=\{0.01,\,1.40,\,4.67,\,7.00\}\,\mathrm{s}$ for $N=2,3,4,5$.}
    \label{fig:convergence115}
\end{figure}

Figures~\ref{fig:convergence85} and \ref{fig:convergence115} show the convergence of the flow map as the polynomial order increases from $N=2$ to $N=5$. In both the expansion and compression cases, the deformation profile clearly converges. The solutions for $N=4$ and $N=5$ are nearly indistinguishable, suggesting that the limiting shape has been reached and validating the convergence of the method.

To visualize the temporal evolution of the flow dynamics and the convergence behaviour discussed above, full animations of the expansion and compression test cases are available at,  \url{https://www.youtube.com/watch?v=r1-CsvZUdQg&t=82s}.

\FloatBarrier

\section{Conclusion}\label{sec:Conclusions}

In this paper, we have presented a novel, higher-order, structure-preserving space-time discretization for the Lagrangian formulation of inviscid barotropic flow. The method is built from first principles, starting with a multisymplectic variational formulation of the continuous action, which is then discretized over a full (2+1D) space-time domain using the tools of Mimetic Spectral Element Methods (MSEM). The central contribution of this work is a formulation that entirely bypasses the mesh distortion and instability issues that plague traditional Lagrangian methods. By defining all kinematic and dynamic variables on a fixed, non-deforming reference configuration ($\mathcal{B}$), all computations are performed on a static mesh. The physical deformation is captured implicitly by the flow map ($\bm{\varphi}$) and its gradient ($\bm{F}$), while forces are mapped back from the physical domain via the first Piola-Kirchhoff stress tensor ($\widetilde{\bm{P}}$).

This structure-preserving design, which leverages a primal-dual staggering of discrete differential forms (cochains), is not merely an approximation; it creates a discrete system that inherits the fundamental conservation laws of the continuous physics. We have shown analytically that the resulting discrete weak formulation exactly conserves discrete linear momentum, angular momentum, and total energy. Numerical experiments for low-Mach number expansion and compression flows validated these theoretical findings. The simulations confirmed the exact conservation of all invariants to machine precision and demonstrated the method's stability and accuracy, with clear convergence of the flow map as the polynomial order was increased.

\newpage

\bibliographystyle{abbrv}
\bibliography{references} 

\end{document}